\documentclass[a4paper,12pt]{amsart}

\usepackage{xcolor}
\usepackage{amsmath}
\usepackage{amsfonts}
\usepackage{amssymb}
\usepackage{amsthm}
\usepackage{amsopn}
\usepackage{amscd}
\usepackage{pstricks}
\usepackage[utf8]{inputenc}
\usepackage{courier}
\usepackage{url}
\usepackage{color}
\usepackage{subfigure}
\usepackage{xypic}
\usepackage{graphicx}
\usepackage{paralist}
\usepackage[all,cmtip]{xy}
\usepackage{tikz}

\usepackage[vcentering,dvips]{geometry}
\entrymodifiers={+!!<0pt,\fontdimen22\textfont2>}
\title{Upgrading and Downgrading Torus Actions}

\author[N.O.~Ilten]{Nathan Owen Ilten}
\address{Max Planck Institut f\"ur Mathematik,
        PF 7280,
        53072 Bonn, Germany}
\email{nilten@cs.uchicago.edu}
\author[R.~Vollmert]{Robert Vollmert}
\address{Institut f\"ur Mathematik und Informatik,
        Freie Universit\"at Berlin,
        Arnimallee 3,
        14195 Berlin, Germany}
\email{vollmert@math.fu-berlin.de}

\definecolor{lcolor}{rgb}{1.0,0.3,0.0}

\newcommand{\p}{\mathfrak{p}}
\newcommand{\CC}{\mathbb{C}}
\newcommand{\QQ}{\mathbb{Q}}
\newcommand{\ZZ}{\mathbb Z}
\newcommand{\NN}{\mathbb N}
\renewcommand{\AA}{\mathbb{A}}
\newcommand{\PP}{\mathbb{P}}

\newcommand{\D}{\mathcal D}
\newcommand{\E}{\mathcal E}
\newcommand{\A}{\mathcal A}
\newcommand{\C}{\mathcal C}
\newcommand{\R}{\mathcal R}
\newcommand{\V}{\mathcal V}
\newcommand{\N}{\mathcal N}
\newcommand{\cX}{\mathcal{X}}
\newcommand{\mcP}{\mathcal P}
\newcommand{\CO}{\mathcal O}
\newcommand{\wt}{\omega}
\newcommand{\dfan}{\mathcal{S}}

\newcommand{\tD}{\widetilde{\mathcal{D}}}
\newcommand{\tT}{\widetilde{T}}
\newcommand{\twt}{\widetilde{\wt}}
\newcommand{\tM}{\widetilde{M}}
\newcommand{\tN}{\widetilde{N}}
\newcommand{\tY}{\widetilde{Y}}

\newcommand{\tsigma}{\widetilde{\sigma}}
\newcommand{\xvert}{\mathbf{vert}}
\newcommand{\xray}{\mathbf{ray}}
\newcommand{\lin}{\mathbf{lin}}
\newcommand{\cox}{\mathbf{cox}}

\newcommand{\hD}{\overline{\mathcal{D}}}
\newcommand{\hT}{\overline{T}}
\newcommand{\hwt}{\overline{\wt}}
\newcommand{\hM}{\overline{M}}
\newcommand{\hN}{\overline{N}}
\newcommand{\hY}{\overline{Y}}
\newcommand{\hu}{\overline{u}}

\newcommand{\sfrac}[1]{{\textstyle \frac{1}{#1}}}
\newcommand{\gquot}{/\!\!/}
\newcommand{\pair}[2]{\langle #1, #2 \rangle}

\DeclareMathOperator{\relint}{relint}
\DeclareMathOperator{\Cl}{Cl}
\DeclareMathOperator{\Pic}{Pic}
\DeclareMathOperator{\NE}{NE}
\DeclareMathOperator{\Spec}{Spec}
\DeclareMathOperator{\Proj}{Proj}
\DeclareMathOperator{\X}{X}
\DeclareMathOperator{\rk}{rk}
\DeclareMathOperator{\WDiv}{WDiv}
\DeclareMathOperator{\CDiv}{CaDiv}
\DeclareMathOperator{\Div}{div}
\DeclareMathOperator{\supp}{supp}

\DeclareMathOperator{\spec}{Spec}
\DeclareMathOperator{\ord}{ord}
\DeclareMathOperator{\conv}{conv}
\DeclareMathOperator{\tail}{tail}

\DeclareMathOperator{\pos}{pos}
\DeclareMathOperator{\slope}{slope}
\DeclareMathOperator{\loc}{loc}
\DeclareMathOperator{\id}{id}

\DeclareMathOperator{\pr}{pr}
\DeclareMathOperator{\pt}{pt}
\DeclareMathOperator{\Cox}{Cox}
\DeclareMathOperator{\Bl}{Bl}

\theoremstyle{theorem}
\newtheorem{theorem}{Theorem}[section]
\newtheorem{cor}{Corollary}[theorem]
\newtheorem{prop}[theorem]{Proposition}
\newtheorem{lemma}[theorem]{Lemma}
\theoremstyle{definition}
\newtheorem{definition}[theorem]{Definition}
\newtheorem{rem}[theorem]{Remark}
\newtheorem{ex}[theorem]{Example}

\newcommand{\minkda}{
\psset{unit=1.4cm}
\begin{pspicture}(-2,-.5)(1,2)
\pspolygon[fillstyle=solid,fillcolor=lightgray,linecolor=lightgray](-.5,.5)(-1.5,1.5)(.5,1.5)(-.5,.5)
\psline{->}(-.5,.5)(-1.6,1.6)
\psline{->}(-.5,.5)(.6,1.6)
\psdot(-.5,.5)
{\scriptsize{
\rput(-.5,.2){$(-\frac{1}{2},\frac{1}{2})$}
}}

\end{pspicture}
}
\newcommand{\minkdb}{
\psset{unit=1.4cm}
\begin{pspicture}(-1.5,-.5)(2.5,2)
\pspolygon[fillstyle=solid,fillcolor=lightgray,linecolor=lightgray](0,0)(1,0)(2,1)(-1,1)(0,0)
\psline{->}(0,0)(-1.1,1.1)
\psline(1,0)
\psline{->}(1,0)(2.1,1.1)
\psdots(0,0)(1,0)
{\scriptsize{
\rput(0,-.3){$(0,0)$}
\rput(1,-.3){$(1,0)$}
}}

\end{pspicture}
}

\newcommand{\minkdecomp}{
\psset{unit=1.4cm}
\begin{pspicture}(-1.5,-.5)(6,1.5)
\pspolygon[fillstyle=solid,fillcolor=lightgray,linecolor=lightgray](0,0)(1.1,1.1)(-1.1,1.1)(0,0)
\psline{->}(0,0)(-1.2,1.2)
\psline{->}(0,0)(1.2,1.2)
\psline[linestyle=dashed,linecolor=gray](-1,.5)(1,.5)
\psline[linewidth=.07,linecolor=gray](-.5,.5)(.5,.5)
\psline[linewidth=.07,linecolor=gray](5,0)(6,0)
\psline[linecolor=gray,linewidth=.07](2.96,.5)(3.04,.5)
\psdots(-1,0)(0,0)(1,0)(-1,1)(0,1)(1,1)
\psdots(2.5,0)(3.5,0)(2.5,1)(3.5,1)
\psdots(5,0)(6,0)(5,1)(6,1)
\rput(1.75,.5){$=$}
\rput(4.25,.5){$+$}
{\scriptsize{
\rput(3,.8){$(-\frac{1}{2},\frac{1}{2})$}
\rput(5,-.3){$(0,0)$}
\rput(6,-.3){$(1,0)$}
\rput(-1.35,.85){$(-1,1)$}
\rput(1.3,.85){$(1,1)$}
}}

\end{pspicture}
}

\newcommand{\fanslicea}{
\psset{unit=1cm}
\begin{pspicture}(-1,-1)(4,3)
\psdots(0,0)(2,0)(2,2)(0,2)(1,1)
\psline(0,0)(2,0)(2,2)(0,2)(0,0)
\psline(2,2)
\psline(0,2)(2,0)
\rput(0,-.3){$v_1$}
\rput(2,-.3){$v_2$}
\rput(0,2.3){$v_3$}
\rput(2,2.3){$v_4$}
\rput(1,.5){$\frac{1}{2}v_5$}
\end{pspicture}
}

\newcommand{\fansliceb}{
\psset{unit=1cm}
\begin{pspicture}(-2,-1)(1,1)
\psdots(0,0)(0,1)(0,2)
\psline(0,0)(0,2)
\psline[doubleline=true,doublesep=.1cm,linestyle=dashed]{->}(-3,1)(-1,1)
\rput(.7,0){\scriptsize{$(1,0)$}}
\rput(.7,1){\scriptsize{$\frac{1}{2}(1,1)$}}
\rput(.7,2){\scriptsize{$(0,1)$}}
\end{pspicture}
}

\newcommand{\nonsnaketail}{
\psset{unit=1cm}
\begin{pspicture}(-3,-1)(3,1)
\psline{<-]}(-2,0)(0,0)
\psline{[->}(2,0)
\rput(-2.5,0){$\rho_2$}
\rput(2.5,0){$\rho_1$}
\rput(0,.4){$0$}
\end{pspicture}
}

\newcommand{\nonsnakeslice}{
\psset{unit=1cm}
\begin{pspicture}(-3,-1)(3,1)
\psline{<-]}(-2,0)(-1,0)
\psline{[->}(-1,0)(2,0)
\rput(-1.2,.4){$-1$}
\end{pspicture}
}

\newcommand{\nonsnakeuptail}{
\psset{unit=.7cm}
\begin{pspicture}(-1.5,-2.5)(3.5,3.5)
\psgrid[gridwidth=1pt,griddots=5,subgriddiv=1,gridlabels=5pt,gridlabelcolor=white](-1,-2)(3,3)
\pspolygon[fillstyle=crosshatch,hatchsep=.1,hatchcolor=gray,linecolor=gray](3,0)(0,0)(1.5,3)(3,3)
\pspolygon[fillstyle=crosshatch,hatchsep=.1,hatchcolor=lightgray,linecolor=lightgray](3,0)(0,0)(0,-2)(3,-2)
\psdot(0,0)
\rput(-.8,0){$(0,0)$}
\end{pspicture}}

\newcommand{\nonsnakeupcoeff}{
\psset{unit=.7cm}
\begin{pspicture}(-1.5,-2.5)(3.5,3.5)
\psgrid[gridwidth=1pt,griddots=5,subgriddiv=1,gridlabels=5pt,gridlabelcolor=white](-1,-2)(3,3)
\pspolygon[fillstyle=crosshatch,hatchsep=.1,hatchcolor=gray,linecolor=gray](3,-1)(0,-1)(2,3)(3,3)
\pspolygon[fillstyle=crosshatch,hatchsep=.1,hatchcolor=lightgray,linecolor=lightgray](3,-1)(0,-1)(0,-2)(3,-2)
\psdot(0,-1)
\rput(-1,-1){$(0,-1)$}
\end{pspicture}}

\newcommand{\grassfantail}{
\psset{unit=.7cm}
\begin{pspicture}(-1.5,-2)(2.5,3)
\psgrid[gridwidth=1pt,griddots=5,subgriddiv=1,gridlabels=5pt,gridlabelcolor=white](-1,-1)(2,2)
\pspolygon[fillstyle=crosshatch,hatchcolor=lightgray,hatchsep=.2,linecolor=lightgray](1,-1)(-1,1)(-1,2)(2,2)(2,-1)
\psline{->}(1,-1)
\psline{->}(2,0)
\psline{->}(2,2)
\psline{->}(0,2)
\psline{->}(-1,1)
\psdots[dotstyle=o](-1,-1)(1,-1)(0,0)(2,0)(-1,1)(1,1)(0,2)(2,2)
\rput(1,-1.4){$\rho_1$}
\rput(2.4,0){$\rho_2$}
\rput(2.3,2.3){$\rho_3$}
\rput(0,2.4){$\rho_4$}
\rput(-1.4,1){$\rho_5$}
\end{pspicture}}

\newcommand{\grassfana}{
\psset{unit=.7cm}
\begin{pspicture}(-2.5,-2)(2.5,4)
\psgrid[gridwidth=1pt,griddots=5,subgriddiv=1,gridlabels=5pt,gridlabelcolor=white](-2,-1)(2,3)
\pspolygon[fillstyle=crosshatch,hatchcolor=lightgray,hatchsep=.2,linecolor=lightgray](-2,2)(-2,3)(2,3)(2,-1)(1,-1)
\psdots[dotstyle=o](-1,-1)(1,-1)(-2,0)(0,0)(2,0)(-1,1)(1,1)(-2,0)(0,2)(2,2)(-1,3)(1,3)
\psline(-2,2)(-1,1)(-1,3)
\psline(1,3)(-1,1)(0,0)(2,2)
\psline(2,0)(0,0)(1,-1)
\psdots(0,0)(-1,1)
\rput(-.2,-.4){$0$}
\end{pspicture}}

\newcommand{\grassfanb}{
\psset{unit=.7cm}
\begin{pspicture}(-2.5,-2)(2.5,4)
\psgrid[gridwidth=1pt,griddots=5,subgriddiv=1,gridlabels=5pt,gridlabelcolor=white](-2,-1)(2,3)
\pspolygon[fillstyle=crosshatch,hatchcolor=lightgray,hatchsep=.2,linecolor=lightgray](-2,2)(-2,3)(2,3)(2,-1)(1,-1)
\psdots[dotstyle=o](-1,-1)(1,-1)(-2,0)(0,0)(2,0)(-1,1)(1,1)(-2,0)(0,2)(2,2)(-1,3)(1,3)
\psline(-2,2)(-1,1)(-1,3)
\psline(1,3)(-1,1)(0,0)(2,2)
\psline(2,0)(0,0)(1,-1)
\psdots(0,0)(-1,1)
\rput(-1.2,.6){$0$}
\end{pspicture}}

\newcommand{\grassfanc}{
\psset{unit=.7cm}
\begin{pspicture}(-1.5,-2)(3.5,4)
\psgrid[gridwidth=1pt,griddots=5,subgriddiv=1,gridlabels=5pt,gridlabelcolor=white](-1,-1)(3,3)
\pspolygon[fillstyle=crosshatch,hatchcolor=lightgray,hatchsep=.2,linecolor=lightgray](1,-1)(-1,1)(-1,3)(3,3)(3,-1)
\psline(1,-1)(-1,1)
\psline(3,3)
\psline(1,3)(1,1)(3,1)
\psdots[dotstyle=o](-1,-1)(1,-1)(3,-1)(0,0)(2,0)(-1,1)(1,1)(3,1)(0,2)(2,2)(-1,3)(1,3)(3,3)
\psdots(0,0)(1,1)
\rput(-.3,-.3){$0$}
\end{pspicture}}

\newcommand{\grasscoeffa}{
\psset{unit=.7cm}
\begin{pspicture}(-1.5,-2)(2.5,3)
\psgrid[gridwidth=1pt,griddots=5,subgriddiv=1,gridlabels=5pt,gridlabelcolor=white](-1,-1)(2,2)
\pspolygon[fillstyle=crosshatch,hatchcolor=lightgray,hatchsep=.2,linecolor=lightgray](0,2)(0,1)(2,1)(2,2)
\psdot(0,1)
\psline{<->}(0,2)(0,1)(2,1)
\rput(0,.7){$e_2$}
\end{pspicture}}

\newcommand{\grasscoeffb}{
\psset{unit=.7cm}
\begin{pspicture}(-1.5,-2)(2.5,3)
\psgrid[gridwidth=1pt,griddots=5,subgriddiv=1,gridlabels=5pt,gridlabelcolor=white](-1,-1)(2,2)
\pspolygon[fillstyle=crosshatch,hatchcolor=lightgray,hatchsep=.2,linecolor=lightgray](-1,2)(-1,0)(2,0)(2,2)
\psdot(-1,0)
\psline{<->}(-1,2)(-1,0)(2,0)
\rput(-1,-.3){$-e_1$}
\end{pspicture}}

\newcommand{\grasscoeffc}{
\psset{unit=.7cm}
\begin{pspicture}(-1.5,-2)(2.5,3)
\psgrid[gridwidth=1pt,griddots=5,subgriddiv=1,gridlabels=5pt,gridlabelcolor=white](-1,-1)(2,2)
\pspolygon[fillstyle=crosshatch,hatchcolor=lightgray,hatchsep=.2,linecolor=lightgray](0,2)(0,-1)(2,-1)(2,2)
\psdot(0,-1)
\psline{<->}(0,2)(0,-1)(2,-1)
\rput(-.6,-1){$-e_2$}
\end{pspicture}}

\newcommand{\grasscoeffd}{
\psset{unit=.7cm}
\begin{pspicture}(-1.5,-2)(2.5,3)
\psgrid[gridwidth=1pt,griddots=5,subgriddiv=1,gridlabels=5pt,gridlabelcolor=white](-1,-1)(2,2)
\pspolygon[fillstyle=crosshatch,hatchcolor=lightgray,hatchsep=.2,linecolor=lightgray](1,2)(1,0)(2,0)(2,2)
\psdot(1,0)
\psline{<->}(1,2)(1,0)(2,0)
\rput(.6,0){$e_1$}
\end{pspicture}}

\newcommand{\grasscoeffe}{
\psset{unit=.7cm}
\begin{pspicture}(-1.5,-2)(2.5,3)
\psgrid[gridwidth=1pt,griddots=5,subgriddiv=1,gridlabels=5pt,gridlabelcolor=white](-1,-1)(2,2)
\pspolygon[fillstyle=crosshatch,hatchcolor=lightgray,hatchsep=.2,linecolor=lightgray](-1,2)(-1,0)(0,-1)(2,-1)(2,2)
\psdots(0,-1)(-1,0)
\psline{<->}(-1,2)(-1,0)(0,-1)(2,-1)
\rput(-1.2,-.3){$-e_1$}
\rput(-.5,-1.2){$-e_2$}
\end{pspicture}}

\newcommand{\grasscube}{
\psset{unit=.35cm}
\begin{pspicture}(0,0)(16,16)
\psdots(0.00,0.00)(10.00,0.00)(0.00,10.00)(10.00,10.00)(6.00,4.00)(16.00,4.00)(6.00,14.00)(16.00,14.00)
\psline(0.00,0.00)(10.00,0.00)(10.00,10.00)(0,10)(0,0)
\psline(10,0)(16.00,4.00)(16.00,14.00)(10,10)
\psline(0,10)(6,14)(16,14)
\psline[linestyle=dashed](0,0)(6,4)(16,4)
\psline[linestyle=dashed](6,4)(6,14)
\psline[linecolor=gray,linewidth=.29](0,0)(10,0)
\psline[linecolor=gray,linewidth=.29](16,14)(10,10)
\psline[linecolor=gray,linewidth=.29,linestyle=dashed](6,4)(6,14)
\rput(5.00,15.00){$e_2$}
\rput(12.00,0.00){$e-e_2$}
\rput(18.00,14.00){$e-e_1$}
\rput(17.00,4.00){$e_4$}
\rput(-2.00,10.00){$e-e_4$}
\rput(9.80,11.300){$e_3$}
\rput(6.50,3.00){$e-e_3$}
\rput(-1.00,0.00){$e_1$}
\rput(5,8){$\color{gray} \Delta_0^c$}
\rput(13.5,10){$\color{gray} \Delta_1^c+e$}
\rput(5,-1){$\color{gray} \Delta_\infty^c$}
\end{pspicture}
}

\newcommand{\tvardefBase}{%
\begin{tikzpicture}[scale=1.5]

\draw[thick] (-0.7,-0.7) -- (0.7,0.7) node[above]{$D_1$};
\draw[thick] (0,-0.7) -- (0,0.7) node[above]{$D_0$};
\draw[thick] (1.3,-0.7) -- (1.3,0.7) node[above]{$D_\infty$};
\draw[dashed] (-0.8,0) -- (1.6,0);

\draw[->] (1.8,0) -- (2.2,0);

\begin{scope}[xshift=2.4cm]
  \draw[thick] (0,-0.7) -- (0,0.7);
  \node at (0,0) [inner sep=0.3mm,circle,fill] {};
  \node[right] at (0,0) {$0$};
  \node at (0,0.5) [inner sep=0.3mm,circle,fill] {};
  \node[right] at (0,0.5) {$s$};
\end{scope}

\end{tikzpicture}
}%

\newcommand{\tvardefDivisors}{%
\begin{tikzpicture}[scale=1.5]

\foreach \y/\lbl in { 0.3/$\mathcal{D}_s$,
                      -0.3/$\mathcal{D}_0$ }
  \draw[thick] (1.6,\y) -- (-0.7,\y) node[left]{\lbl};

\foreach \x/\lbl in {
      0/$\Delta_0$/,
      0.5/$\Delta_1$/,
      1.3/$\Delta^-$/
    }
{
  \node[inner sep=0.3mm,circle,fill] at (\x,0.3) {};
  \node[above] at (\x,0.3) {\lbl};
}

\foreach \x/\lbl in {
      0/$\Delta_0 + \Delta_1$/,
      1.3/$\Delta^-$/
    }
{
  \node[inner sep=0.3mm,circle,fill] at (\x,-0.3) {};
  \node[below] at (\x,-0.3) {\lbl};
}

\foreach \x/\lbl in {0/$0$, 0.5/$s$/,
                     1.3/$\infty$}
  \node at (\x,0) {\lbl};

\end{tikzpicture}
}%

\begin{document}
\maketitle
\begin{abstract}
K. Altmann and J. Hausen have shown that affine $T$-varieties can be described in terms of p-divisors. Given a p-divisor describing a $T$-variety $X$, we show how to construct new p-divisors describing $X$ with respect to actions by larger tori. Conversely, if $\dim T=\dim X-1$, we show how to construct new p-divisors describing $X$ with respect to actions by closed subtori of $T$. As a first application, we give explicit constructions for the p-divisors describing certain Cox rings. Furthermore, we show how to upgrade the p-divisors describing the total spaces of homogeneous deformations of toric varieties.
\end{abstract}

\section*{Introduction}
Recall that a normal variety $X$ together with an effective action by some algebraic torus $T$ is called a \emph{$T$-variety}. The \emph{complexity} of a $T$-variety $X$ is equal to $\dim X-\dim T$. Complexity zero $T$-varieties are thus just the well-studied toric varieties.

In \cite{altmann:06a}, K. Altmann and J. Hausen introduced an approach to studying affine $T$-varieties which generalizes the correspondence between affine toric varieties and polyhedral cones. Given some affine $T$-variety $X$ of complexity $k$, the basic idea is to split the information describing $X$ into a $k$-dimensional variety $Y$ together with a divisor $\D$ on $Y$ having polyhedral instead of integral coefficients.  The variety $Y$ is a sort of quotient of $X$ by $T$. The divisor $\D$ encodes the information lost by passing from $X$ to $Y$ and is called a \emph{p-divisor}, see Definition \ref{def:pdiv}. In particular, any p-divisor $\D$ determines a $T$-variety which we denote by $X(\D)$, see Theorem \ref{thm:pdivmain}. This basic idea globalizes to non-affine $T$-varieties as well, see \cite{altmann:08a}.  For a survey of the many results to which this approach leads, see \cite{tsurvey}.

The primary goal of this paper is to describe p-divisors arising when
upgrading the action on an affine $T$-variety to an action by a larger torus.
  Here, we have some $T$-variety
$X=X(\D)$ for some p-divisor $\D$ on a normal variety
$Y$. Suppose that $Y$ itself is a
$T'$-variety, and $\D$ is a $T'$-invariant p-divisor; we make this latter notion precise in Section \ref{sec:upgrade}. Then $X$ is in fact a $\tT$-variety where
$\tT:=T\times T'$. In this situation, we wish to find a
new p-divisor $\tD$ such that $X=X(\tD)$ as a $\tT$ variety.
After some mild assumptions on $Y$, we are able to do this completely in combinatorial terms, see Section \ref{sec:upgrade} and Theorem \ref{thm:tD}. The proof of Theorem \ref{thm:tD} requires an analysis of semiample divisors on $T$-varieties of arbitrary complexity, which we carry out in Section \ref{sec:ssa}.

We are also interested in a sort of inverse problem, namely p-divisors arising via downgrading, that is, considering the action of some subtorus.
 Here, we have some $T$-variety
$X=X(\D)$ for some p-divisor $\D$ on a semiprojective normal variety
$Y$ along with a subtorus $\hT\subset T$. In this situation, we wish to find some new p-divisor $\hD$ such that $X=X(\hD)$ as a $\hT$-variety.
This downgrading procedure was described explicitly in Section 11 of \cite{altmann:06a} for the case that $\dim T=\dim X$, that is, $X$ is a toric variety. In Section \ref{sec:downgradeone} we present a similar description for the complexity-one case, that is, when $\dim T=\dim X-1$.
A key ingredient is the use of \emph{divisorial polyhedra} (see Section \ref{sec:divpoly}), a generalization of the divisorial polytopes found in \cite{ilten:09d}.
We also discuss problems arising when downgrading actions of higher complexity, see Section \ref{sec:downgradetwo}.

We then turn to applications of our upgrading procedure. The spectrum $X^{\cox}$ of the Cox ring of some Mori dream space $X$ has a natural effective action by the so-called Picard torus $T$, see Section \ref{sec:cox}. If $X$ itself is a $T'$-variety, then the big torus $\tT=T\times T'$ also acts on $X^\cox$. We use our upgrade procedure to construct a p-divisor describing $X^{\cox}$ as a $\tT$-variety in the case that the divisor class group of $X$ is torsion-free, and the Chow quotient of $X$ has Picard group equal to $\ZZ$. In particular, we partially recover the result of K. Altmann and L. Petersen in \cite{altmann:10a}, which describes the
p-divisor for the Cox ring of a rational complexity-one $T$-variety.

 In a second application, we consider the  homogeneous deformations of affine complexity-one $T$-varieties constructed by the present authors in \cite{ilten:09b}.
If the special fiber of such a deformation is toric, then this bigger torus acts on the total space of the deformation as well. We again apply our upgrade procedure to explicitly calculate the corresponding p-divisor in Section \ref{sec:def}.
\\
\\
\noindent{\bf Acknowledgements:} The authors would like to thank Piotr Achinger, Klaus Altmann, and Hendrik S\"u\ss{} for helpful discussions, as well as the MPIM in Bonn for support and a pleasant working environment.

\section{$T$-Varieties and Invariant Divisors}\label{sec:divisors}
We begin by recalling some basics on polyhedral divisors and affine $T$-varieties which can be found in \cite{altmann:06a}.
Let $N$ be a lattice, $M$ its dual, and let $N_\QQ$ and $M_\QQ$ be the associated $\QQ$-vector spaces. Recall that for any polyhedron $\Delta\subset N_\QQ$, its \emph{tailcone} is the set
$$\tail(\Delta):=\{v\in N_\QQ\ |\ v+\Delta\subseteq \Delta\}.$$
Let $Y$ be a normal variety; we will always be working over an algebraically closed characteristic zero field $\CC$. Fix some cone $\sigma\subseteq N_\QQ$. A \emph{polyhedral divisor} on $Y$ with tailcone $\sigma$ is a formal sum
$$
\D=\sum \D_P\otimes P
$$
over all prime divisors $P\subset Y$, where the coefficients $\D_P$ are either polyhedra in $N_\QQ$ with tailcone $\sigma$ or are the empty set, and only finitely many $\D_P$ differ from $\sigma$.

Let $\wt=\sigma^\vee$. Then a polyhedral divisor $\D$ determines an evaluation map $$\D\colon\wt\to\WDiv_{\QQ\cup\{\infty\}} (Y),$$ where $\WDiv_{\QQ\cup \{\infty\}}(Y)$ are Weil divisors on $Y$ with coefficients in $\QQ\cup\{\infty\}$. Indeed, we set
$$
\D(u):=\sum_{P\in Y} \inf \pair{\D_P}{u} \cdot P
$$
for any $u\in\wt$, where the infimum over the empty set is $\infty$. This map is piecewise-linear and \emph{convex}, i.e. $\D(u)+\D(u')\leq \D(u+u')$ for $u,u'\in \wt$. In fact, there is a bijection between polyhedral divisors on $Y$ with tailcone $\sigma$ and convex piecewise-linear maps $\wt\to\WDiv_{\QQ\cup\{\infty\}}(Y)$, cf. \cite[Proposition 1.5]{altmann:06a}. Thus, we will switch interchangeably between these two descriptions.
While the coefficients of elements of $\WDiv_{\QQ\cup \{\infty\}}(Y)$ might seem slightly odd, one can still define 
$$
L(D):=\{f\in\CC(X)^*\ |\ \Div(f)+D\geq 0\}\cup\{0\}
$$
for any $D\in\WDiv_{\QQ\cup \{\infty\}}(Y)$.
We could alternatively view $D$ as a $\QQ$-divisor on its locus, denoted by $\loc (D)$, which is defined to be the complement of all those prime divisors $P$ on $Y$ for which $D$ has coefficient $\infty$. Note that any time we talk about some property of a divisor $D\in\WDiv_{\QQ\cup\{\infty\}}(Y)$, we are actually referring to the same property of the restriction of $D$ to $\loc (D)$. If $\D$ is a polyhedral divisor, then its locus is defined to be $\loc(\D(u))$ for any $u\in \wt$.

\begin{definition}\label{def:pdiv} 
Let $\D$ be a polyhedral divisor on $Y$ with tailcone $\sigma$.
\begin{enumerate}
	\item $\D$ is \emph{$\QQ$-Cartier} if for all $u\in\wt$, $\D(u)$ is $\QQ$-Cartier;
\item $\D$ is \emph{semiample} if for all $u\in \wt$, $\D$ is \emph{semiample}, that is, a multiple of $\D(u)$ is base-point free;
	\item $\D$ is \emph{big} if $\D(u)$ is \emph{big} for all $u \in \relint \wt$, that is, a multiple of $\D(u)$ admits a section with affine complement;
\item $\D$ is \emph{proper} or simply a \emph{p-divisor} if $\loc (\D)$ is semiprojective,\footnote{Recall that a variety is semiprojective if it is projective over an affine variety.} $\wt$ is full-dimensional, and $\D$ is $\QQ$-Cartier, semiample, and big.
\item If $\D$ is a p-divisor, we call $\wt$ its \emph{weight cone}.\end{enumerate}
\end{definition}

To any polyhedral divisor $\D$ we can associate an affine scheme
$$
X(\D):=\spec\left(\bigoplus_{u\in \wt\cap M} L(\D(u))\cdot\chi^u\right).
$$
This scheme comes with an action by the torus $T=\spec \CC[M]$ and a rational map $X(\D)\dashrightarrow \loc(\D)$.
\begin{theorem}[{\cite[cf. Theorems 3.1 and 3.4]{altmann:06a}}]\label{thm:pdivmain}
If $\D$ is a p-divisor, then $X(\D)$ is a $T$-variety of complexity equal to the dimension of $Y$. Conversely, for any complexity $k$ $T$-variety $X$, there is a p-divisor $\D$ on some normal $k$-dimensional variety with $X(\D)$ equivariantly isomorphic to $X$.
\end{theorem}

To deal with the situation of non-affine $T$-varieties, one must replace the single p-divisor $\D$ by a so-called \emph{divisorial fan} $\dfan$ consisting of a finite set of p-divisors on some common variety $Y$. The idea is that the varieties $X(\D)$ for $\D\in\dfan$ give an invariant affine covering for the $T$-variety of interest. To ensure that these affine varieties glue together in a coherent manner, the elements of $\dfan$ must satisfy certain compatibility conditions. In particular,
for any prime divisor $P$ on $Y$, the set of polyhedra
$$
\dfan_P:=\{\D_P\subset N_\QQ\ | \D\in\dfan\}
$$
should form a polyhedral complex in $N_\QQ$ which we call a slice.
This is the one important property of divisorial fans which we shall need; for a precise definition, the interested reader may refer to \cite[Definition 5.2]{altmann:08a}. In any case, such a set of p-divisors determines a $T$-scheme which we denote by $X(\dfan)$ \cite[Theorem 5.3]{altmann:08a}. 
The tailcones of the polyhedra in any slice $\dfan_P$ form a fan $\tail(\dfan_P)$; this tailfan actually doesn't depend on $P$, so we also denote it simply by $\tail(\dfan)$.
We denote by  $\loc (\dfan)$ the complement in $Y$ of all prime divisors $P$ with $\dfan_P=\emptyset$.

Fix now some  $T$-variety $X$ described by a divisorial fan $\dfan$ on
$Y$ and let $n=\dim T$.
In \cite{petersen:08a}, L. Petersen and H. S\"u\ss{} describe 
all $T$-invariant prime Weil divisors on a $T$-variety. 
Indeed, there are ``vertical'' invariant prime divisors arising as the closure of a family 
of $n$-dimensional $T$-orbits. 
Such divisors are parametrized by prime divisors 
$P\subset Y$ together with $v$ a vertex of $\dfan_P$ satisfying some 
additional condition \cite[Proposition 3.13]{petersen:08a}; let $\xvert_P(\dfan)$ denote the set of such vertices. We denote the divisor corresponding to such 
$P$ and $v$ by $D_{P,v}$. All other invariant prime divisors are ``horizontal'' and arise
as the closure of a family of $(n-1)$-dimensional $T$-orbits. These are parametrized by rays $\rho$ of $\tail(\dfan)$ satisfying an
additional condition; let $\xray(\dfan)$ denote the set of such rays. 
We denote the divisor corresponding to a ray $\rho$ by $D_\rho$.
For any ray $\rho$, we denote its primitive lattice generator by $v_\rho$.

\begin{lemma}\label{lemma:hauptdivisor}[cf. {\cite[Proposition 3.14]{petersen:08a}}]
Consider some $f\in\CC(Y)$, $u\in M$. Then 
$$
\Div(f\cdot \chi^u)=\sum_\rho \pair{v_\rho}{u} D_\rho+\sum_{P,v}\mu(v)(\pair{v}{u} + \ord_P(f))D_{P,v}
$$
where $\mu(v)$ is the smallest positive integer such that $\mu(v)v\in N$.
\end{lemma}

Now consider any $T$-invariant $\QQ$-Weil divisor $D$ on $X$, which by the above description we can write as 
$$
D=\sum_\rho a_\rho D_\rho+\sum_{P,v}\mu(v)b_{P,v}D_{P,v}.
$$
To any such divisor $D$, we associate a polyhedron $\Box^D\subset M_\QQ$ and a piecewise-affine concave function $\Psi^D\colon\Box^D\to\WDiv_{\QQ\cup\{\infty\}} Y$ as follows:

\begin{align*}
	\Box^D:&=\{u\in M_\QQ\ |\ \pair{v_\rho}{u} + a_\rho\geq
        0\ \textrm{for all}\ \rho\in\xray(\dfan)\};\\
	\Psi^D_P(u):&= \inf_{v\in\xvert_P(\dfan)} (\pair{v}{u}+b_{P,v})\in\QQ\cup\{\infty\};\\
\Psi^D(u):&=\sum_P \Psi^D_P(u)P.
\end{align*}
This generalizes the function $h^*$ constructed in \cite{petersen:08a}.
All $\Box^D$ share the same tailcone, dual to the cone generated by the rays $\rho\in\xray(\dfan)$.

\begin{prop}\label{prop:gs}
For $D$ as above, we have
$$
L(D)=\bigoplus_{u\in\Box_D\cap M} L(\Psi^D(u))\cdot\chi^u.
$$
\end{prop}
\begin{proof}
The proof is similar to that of \cite{petersen:08a}, Proposition~3.23. Consider some $u \in M$ and $f\in\CC(Y)$. Then by Lemma
\ref{lemma:hauptdivisor},
$f\cdot\chi^u\in L(D)$ if and only if 
$$
\sum_\rho \pair{v_\rho}{u} D_\rho+\sum_{P,v}\mu(v)(\pair{v}{u} + \ord_P(f))D_{P,v}+D\geq 0.
$$
This is equivalent to satisfying the following inequalities:
\begin{align*}
\pair{v_\rho}{u} +a_\rho&\geq 0;\\
\pair{v}{u} + \ord_P(f)+b_{P,v}&\geq 0.
\end{align*}
The first line of inequalities is equivalent to $u\in\Box^D$. The
second is equivalent to $f\in L(\Psi^D(u))$.
\end{proof}

\begin{ex}
Let $Y$ be the blowup of $\AA^2$ in the origin, and consider the p-divisor $\D\colon\QQ_{\geq 0}\to 0$ on $Y$. Then $X:=X(\D)=\AA^3$. Now, for any invariant divisor $D$ on $X$, $\Psi_E^D=\infty$, where $E$ is the exceptional divisor in $Y$, since one can check that $\xvert_E(\{\D\})=\emptyset$. If we use the above proposition to calculate $L(D)$, we end up calculating global sections of $\QQ$-divisors on $Y\setminus E=\AA^2\setminus \{0\}$, which is not semiprojective. To avoid calculating global sections on such unpleasant varieties, one must make further assumptions about the relevant $T$-variety $X$ and its divisorial fan $\dfan$, see Proposition~\ref{prop:semiproj}.
\end{ex}

We now recall some general facts about $T$-invariant divisors we shall
need later.
\begin{rem}\label{rem:invgen}
Suppose that $D$ is any $T$-invariant $\QQ$-Cartier divisor on $X$. If $\CO(D)$
is globally generated, these generators can be taken to be $T$-invariant.
Indeed, $H^0(X,\CO(D))$ is generated as an $H^0(X,\CO_X)$-module by
$T$-invariant sections $s_1,\ldots,s_k$, which
will then globally generate $\CO(D)$.
\end{rem}

\begin{lemma}[cf. {\cite[page 61]{fulton:93a}}]\label{lemma:invcover}
Suppose that $X=X(\D)$ is an affine $T$-variety and consider any $T$-invariant Cartier divisor $D$ on
$X$. Then there is a $T$-invariant covering of $X$ on which $D$ is
principal and defined by invariant functions.
\end{lemma}
\begin{proof}[Proof (communicated by H. S\"u\ss{})]
It is sufficient to consider the case $D$ effective. Thus, $D$
corresponds to an ideal $I$ of $A:=H^0(X,\CO_X)$ which is
$M$-homogeneous since $D$ is $T$-invariant. Let $f_1,\ldots,f_k$ be
homogeneous generators of $I$. 
Consider some prime $\p\in\spec A$. Then $I_\p\subset A_\p$ is
generated by some $f_j$. Indeed, some $f_j$ doesn't lie in $\p\cdot I$,
otherwise the $f_i$ can't generate $I$. Since $I_\p$ is principal, by
Nakayama's lemma this $f_j$ then generates $I_\p$.

It follows that the $f_i$ locally define $D$, say on
some open sets $U_i$. Now, we can even find a $T$-invariant cover on which
the $f_i$
define $D$.  Indeed, let $U_i'$ be the complement of all prime
divisors where $f_i$ doesn't define $D$. Since $f_i$ and $D$ are
$T$-invariant, then $U_i'$ is as well.  Furthermore, the $U_i'$ cover
$X$ since $U_i\subset U_i'$.
\end{proof}

In general, the study of divisors on $T$-varieties of complexity higher
than one is complicated by the descriptions of $\xray(\dfan)$ and $\xvert_P(\dfan)$.
To simplify things, one may consider a special situation.

\begin{definition}
A divisorial fan $\dfan$ on $Y$ is \emph{contraction-free} if for all $\D \in\dfan$, the locus of $\D$ is affine.
\end{definition}

Equivalently, $\dfan$ is contraction-free if and only if the rational
quotient map $\pi\colon X(\dfan)\dashrightarrow Y$ is regular or,
using the notation of \cite{altmann:06a}, the affine contraction maps
$\widetilde{X}(\D) \to X(\D)$ are isomorphisms.

\begin{prop}\label{prop:contraction-free}
Let $\dfan$ be a divisorial fan on $Y$.  Then there is a contraction-free divisorial
fan $\dfan'$ on $Y$ such that $\dfan_P=\dfan_P'$ for all slices and there is a birational proper map $X(\dfan')\to X(\dfan)$. 
\end{prop}
\begin{proof}
For $\D$ in $\dfan$, let $\dfan^\D$ be a divisorial fan containing
p-divisors with affine loci such that $X(\dfan^\D)=\spec_Y\A(\D)$. One
easily checks that the union of all such divisorial fans induces a
divisorial fan $\dfan'$ with the desired properties, see for example
\cite[Theorem 3.1]{altmann:06a} and \cite[Remark 3.5]{altmann:08a}. 
\end{proof}

\begin{rem}
If $\dfan$ is a contraction-free divisorial fan, then
$\xray(\dfan)$ is the set of all rays in $\tail(\dfan)$, and
$\xvert_P(\dfan)$ is the set of all vertices in $\dfan_P$.  Indeed,
the bigness conditions in \cite[Proposition 3.13]{petersen:08a} are
trivially fulfilled.
This greatly simplifies the study of divisors on contraction-free $T$-varieties.
\end{rem}

\section{Upgrading a Torus Action}\label{sec:upgrade}
 We now consider the situation where we upgrade the torus action on
 some affine $T$-variety. 
 Assume we have some $T$-variety
$X=X(\D)$ for some polyhedral divisor $\D$ on a  normal variety
$Y$ with weight cone $\wt$. Suppose that $Y$ itself is a
$T'$-variety, and $\D(u)$ is a $T'$-invariant divisor for all $u\in
\wt$. Then $X$ is in fact a $\tT$-variety where
$\tT:=T\times T'$. In this situation, we wish to find a
new polyhedral divisor $\tD$ such that $X=X(\tD)$ as a $\tT$ variety.

Let $M$ and $M'$ be the character lattices of respectively $T$ and
$T'$. Since $Y$ is a $T'$-variety, we have $Y=X(\dfan)$ for some
divisorial fan $\dfan$ on a variety $Y'$. 
We then define
$$
\twt:=\{(u,u')\in M_\QQ\oplus M'_\QQ\ | u\in\wt,\ u'\in\Box^{\D(u)}\}
$$ 
and take
\begin{align*}
\tD\colon\twt\to\WDiv_{\QQ\cup\{\infty\}} Y'\\
(u,u')\mapsto\Psi^{\D(u)}(u').
\end{align*}

It will be useful in proving convexity properties to have a dual
view of this construction.
$\D$ can be expressed as a divisor with polyhedral coefficients
$$
\D = \sum_{\rho\in\xray(\dfan)} \Delta_\rho \otimes D_\rho + \sum_{\substack{P\subset Y'\\v\in\xvert_P(\dfan)}} \Delta_{P,v} \otimes {\mu(v)}D_{P,v},
$$
where the coefficients $\Delta$ have the common tailcone $\sigma = \wt^\vee$.

\begin{prop}\label{prop:dualupgrade}
The dual of the upgraded weight cone $\twt$ is
$$
\tsigma = \pos \big\{ (\sigma \times \{0\}) \cup 
                      \bigcup_{\rho\in\xray(\dfan)} (\Delta_\rho \times \{v_\rho\}) \big\},
$$
where $v_\rho$ is the primitive lattice generator of $\rho$. The upgraded divisor
$\tD$ can be expressed dually as $\tD = \sum \Delta_P \otimes P$ with
$$
\Delta_P = \conv\big\{\Delta_{P,v} \times \{v\} \mid v \in \xvert_P(\dfan) \big\} + \tsigma.
$$
\end{prop}

\begin{proof}
For $(u,u')$ to be an element of $\twt$, first $u$ must be in $\wt$,
equivalently $\pair{(v,0)}{(u,u')} \ge 0$ for all $v \in \sigma$.
Then, we need $u' \in \Box^{\D(u)}$, equivalently $\pair{v_\rho}{u'} + a_\rho \ge 0$ for all occurring rays $\rho$, with $a_\rho$ the coefficient of $D_\rho$ in $\D(u)$.
Thus, $a_\rho = \min \pair{\Delta_\rho}{u}$, and the  condition becomes
$$
    \pair{v_\rho}{u'} + \min_{x \in \Delta_\rho} \pair{x}{u}
=   \min_{x \in \Delta_\rho} \pair{(x,v_\rho)}{(u,u')}
\ge 0.
$$
This proves the claim about $\tsigma$.

For the second claim, consider a prime divisor $P$ on $Y'$.
By definition of $\tD$, the coefficient of $P$ in $\tD(u,u')$ can
be calculated as follows.
\begin{align*}
\Psi^{\D(u)}_P(u') &= \min_{v \in \xvert_P(\dfan)} \pair{v}{u'}
                      + b_{P,v} \\
                   &= \min_{v \in \xvert_P(\dfan)} \pair{v}{u'}
                     + \min_{x \in \Delta_{P,v}} \pair{x}{u} \\
                   &= \min_{v \in \xvert_P(\dfan)} \min_{x \in \Delta_{P,v}}
                       \pair{(x,v)}{(u,u')} \\
                   &= \min_{(x,v) \in \Delta_P} \pair{(x,v)}{(u,u')} \rangle
\end{align*}
This is just the evaluation of the polyhedral divisor
$\sum \Delta_P \otimes P$ at $(u,u')$.
\end{proof}

\begin{theorem}\label{thm:tD}
With the above notation, $\twt$ is a polyhedral cone, $\tD$ is a
polyhedral divisor, and $X=X(\tD)$ as a $\tT$-variety. If $\D$ is a p-divisor and $\dfan$ is contraction-free on smooth $Y'$, then $\tD$ is a p-divisor, and $\twt$ is the weight cone for the $\tT$-action on $X$.
\end{theorem}

\begin{rem}
The hypothesis that $\dfan$ be contraction-free on smooth $Y'$ doesn't pose a large problem. Indeed, we can first resolve the singularities of $Y'$ and pull back $\dfan$. By Proposition~\ref{prop:contraction-free} we can then replace the pullback of $\dfan$ by a contraction-free $\dfan'$, and pull back $\D$ to $\D'$ on $X(\dfan')$.
\end{rem}

\begin{rem}
Theorem \ref{thm:tD} can also be used to upgrade non-affine $T$-varieties. Indeed, if $\dfan$ is a contraction-free divisorial fan on smooth $Y$, and $\Xi$ is a divisorial fan on $X(\dfan)$ consisting of invariant p-divisors, then the upgraded p-divisors $\{\tD\}_{\D\in\Xi}$ form a divisorial fan describing $X(\Xi)$ with the upgraded torus action.
\end{rem}

We delay the proof of the theorem to Section \ref{sec:upgradesrev}.
For now, we present a number of examples. Further applications are discussed in Sections~\ref{sec:cox} and~\ref{sec:def}.

\begin{ex}[Affine cones]
Let $X=X(\dfan)$ be a $T'$-variety, with $\dfan$ a contraction-free divisorial fan on some smooth $Y$, and let $D=\sum a_\rho D_\rho+\sum b_{P,v}\mu(v)D_{P.v}$ be some very ample invariant Cartier divisor on $X$ giving a projectively normal embedding $X\subset \PP^m$. Then the cone $C(X)\subset\AA^{m+1}$ over $X$ is a $\CC^*$-variety given by the p-divisor $\D=[1,\infty)\otimes D$. By the above theorem, $C(X)$ is also a $\tT=\CC^*\times T'$-variety given by the p-divisor  
$$
\tD=\sum_P (\conv\{\{b_{P,v}\}\times\{v\}\}+\tsigma)\otimes P
$$
where $\tsigma=\QQ_{\geq0}\times\{0\}+\pos \{\{a_\rho\}\times \{v_\rho\}\}$. This generalizes Proposition~4.1 of \cite{ilten:09d}.
\end{ex}

\begin{figure}
\begin{center}
\subfigure[$\tail(\dfan)$]{\nonsnaketail}
\subfigure[$\dfan_\infty$]{\nonsnakeslice}
\end{center}
\caption{A divisorial fan for $\PP^2$}\label{fig:nonsnakefan}
\end{figure}

\begin{figure}
\begin{center}
\subfigure[$\tsigma$]{\nonsnakeuptail}
\subfigure[$\Delta_\infty$]{\nonsnakeupcoeff}
\end{center}
\caption{Upgrading with non-contraction-free $\dfan$}\label{fig:nonsnakeup}
\end{figure}

\begin{ex}[A non-contraction-free $\PP^2$]\label{ex:nonsnake}
Consider some divisorial fan $\dfan$ on $Y'=\PP^1$ with tailfan and single nontrivial slice $\dfan_\infty$ as pictured in Figure \ref{fig:nonsnakefan}, where all polyhedra with tailcone $\rho_2$ belong to the same polyhedral divisor, but those with tailcone $\rho_1$ don't.  The resulting $T'$-variety $Y=X(\dfan)$ is in fact $\PP^2$, but $\dfan$ is not contraction-free. Here, $\xray(\dfan)$ consists of $\rho_1$. Now consider the p-divisor $\D=\Delta_{\rho_1}\otimes D_{\rho_1}$, where $\Delta_{\rho_1}=[\sfrac{2},\infty)\subset N_\QQ=\ZZ_\QQ$.

In the upgraded lattice $\tN=N\oplus N'=\ZZ^2$, the upgraded tailcone $\tsigma$ is generated by $(1,0)$ and $(1,2)$, see  the darkly shaded region of Figure \ref{fig:nonsnakeup}(a). Likewise, the upgraded polyhedral divisor $\tD$ on $\PP^1$ is $\Delta_\infty\otimes \{\infty\}$, where $\Delta_\infty=(0,-1)+\tsigma$, see the darkly shaded region of Figure \ref{fig:nonsnakeup}. $\tD$ is clearly not proper.

However, we can replace $\dfan$ with the contraction-free $\dfan'$ gotten by requiring that not all polyhedra with tailcone $\rho_2$ belong to the same polyhedral divisor. Now, $\xray(\dfan)$ consists of $\rho1$ and $\rho_2$. The pullback of $\D$ to $X(\dfan')$ is still $\Delta_{\rho_1}\otimes D_{\rho_1}$, but the upgraded tailcone and coefficient now encompass the lightly and darkly shaded regions of Figure \ref{fig:nonsnakeup}. In particular, the resulting polyhedral divisor is a p-divisor. 
\end{ex}

\begin{figure}
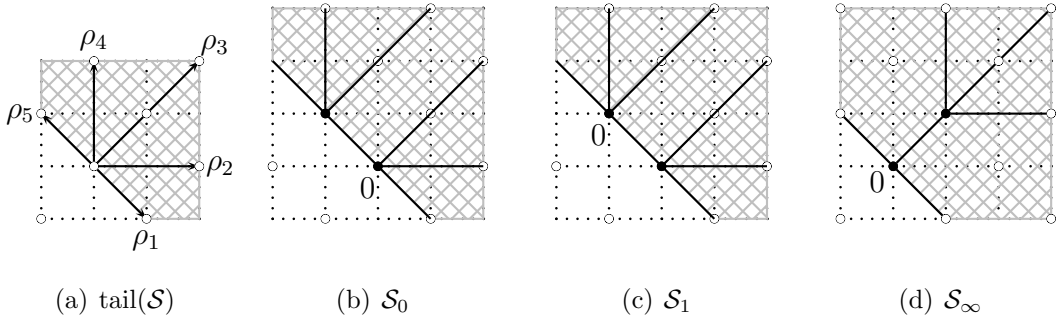

\subfigure[$\tail(\dfan)$]{\grassfantail}
\subfigure[$\dfan_0$]{\grassfana}
\subfigure[$\dfan_1$]{\grassfanb}
\subfigure[$\dfan_\infty$]{\grassfanc}
\caption{A contraction-free divisorial fan on $\PP^1$}\label{fig:grassfan}
\end{figure}

\begin{figure}
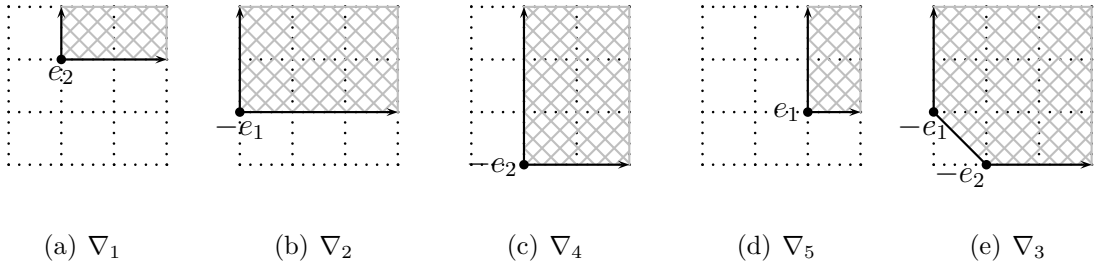

\subfigure[$\nabla_1$]{\grasscoeffa}
\subfigure[$\nabla_2$]{\grasscoeffb}
\subfigure[$\nabla_4$]{\grasscoeffc}
\subfigure[$\nabla_5$]{\grasscoeffd}
\subfigure[$\nabla_3$]{\grasscoeffe}
\caption{Coefficients for a p-divisor}\label{fig:grasscoeff}
\end{figure}

\begin{figure}
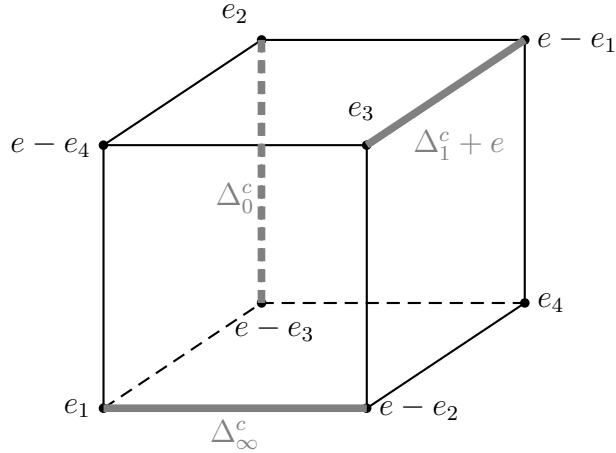

\begin{center}
\grasscube
\end{center} 
\caption{A p-divisor for the cone over $G(2,4)$}\label{fig:cube}
\end{figure}

\begin{ex}[The cone over $G(2,4)$]\label{ex:grassupgrade}
Let $Y'=\PP^1$, and let $N'\subset\ZZ^2$ be the rank two sublattice of $\ZZ^2$ generated by $2f_1,2f_2,f_1+f_2$, where $f_1,f_2$ is the standard basis of $\ZZ^2$. We then consider any contraction-free divisorial fan $\dfan$ in $N_\QQ'$ with slices as pictured in Figure \ref{fig:grassfan}, where dots  denote lattice points of $N'$. Now let $\D$ be the polyhedral divisor on $Y=X(\dfan)$ given by
\begin{align*}\D=\nabla_1\otimes(D_{\rho_1}+D_{0,0})+\nabla_2\otimes(D_{\rho_2}+D_{1,0}+D_{1,f_1-f_2})
+\nabla_4\otimes(D_{\rho_4}+D_{\infty,f_1+f_2})\\
+\nabla_5\otimes(D_{\rho_5}+D_{0,f_2-f_1}+D_{\infty,0})+\nabla_3\otimes D_{\rho_3}
\end{align*}
where $\nabla_i\subset N_\QQ$ are the polyhedra pictured in Figure \ref{fig:grasscoeff}; let $e_1,e_2$ be the standard basis for $N=\ZZ^2$. An easy calculation shows that $\D$ is in fact a p-divisor.

We can thus consider $X=X(\D)$, and try to upgrade the action of $T^N$ on $X$ to that of $\tT=T^{N\oplus N'}$. Let $\tN$ be the sublattice of $\QQ^4$ generated by the standard basis $e_1,e_2,e_3,e_4$ together with $e=\sfrac{2}\sum e_i$. Then $N\oplus N'\cong \tN$, with an isomorphism given by
\begin{align*}
\phi\colon N\oplus N'&\to \tN\\
e_i&\mapsto e_i\\
f_i&\mapsto \frac{1}{2}(e_i+e_{i+2}). 
\end{align*}
Viewing the upgraded tailcone and coefficients of $\tD$ in $\tN$, we have that $\tsigma$ is the cone over the cube of Figure \ref{fig:cube}. Likewise, the upgraded p-divisor is 
$$
\tD=(\Delta_0^c+\tsigma)\otimes\{0\}+
(\Delta_1^c+\tsigma)\otimes\{1\}+
(\Delta_\infty^c+\tsigma)\otimes\{\infty\}
$$
where the $\Delta_P^c$ are the gray line segments pictured in Figure \ref{fig:cube}.  From this description, we recognize that our $X=X(\D)=X(\tD)$ is the cone over $G(2,4)$ in its Pl\"ucker embedding, see \cite[Example 10.6]{altmann:03a} or \cite[Section 6]{altmann:08b}.
\end{ex}

\begin{rem}\label{rem:converse}
Let $X$ be a $\tT$-variety, $T\subset\tT$ some closed subtorus, and $T'=\tT/T$ the quotient torus. Suppose that $\D$ is some p-divisor on a variety $Y$ with $X(\D)=X$ as a $T$-variety.  We know that we should be able to upgrade the action of $T$ on $X$ to that of $\tT$; the question is whether or not we are able to apply Theorem \ref{thm:tD} to our $\D$ and $Y$ to get an upgraded p-divisor $\tD$.

In general, the answer is negative. In fact, $Y$ itself need not admit a $T'$-action. To see this, suppose that $Y$ does in fact admit a $T'$-action. Then we can blow up $Y$ in any point which is not fixed by $T'$, which will destroy the $T'$-action.

However, in the above situation, one can always find a $Y$ with $T'$-action and $T'$-invariant p-divisor $\D$. Indeed, we can take $Y$ to be the normalization of the Chow quotient and $\D$ a p-divisor as in the proof of \cite[Theorem 3.4]{altmann:06a}. We leave it to the reader to check that this construction preserves the residual $T'$-action.
\end{rem}

\begin{ex}[A $\CC^3$-like threefold]
Consider the hypersurface $X$ given by the equation $x+x^2y+z^2+t^3$ in $\AA^4$. $X$ admits a $\CC^*$-action induced by the action on $\AA^4$ given by weights $(6,-6,3,2)$ on the coordinates $(x,y,z,t)$. $X$ is very close to being $\CC^3$, for example, it is smooth, factorial, and birational and diffeomorphic to $\CC^3$. However, the $\CC^*$-action cannot be upgraded to an effective action by a two-dimensional torus. See for example \cite{makar:96a} for a discussion of this hypersurface.

We can also use Remark \ref{rem:converse} to  see that $X$ admits at most a complexity-two action. Section 11 of \cite{altmann:06a} tells us how to calculate a p-divisor  for $X$ with respect to the above $\CC^*$-action. Indeed, let $Y$ be the blowup of $\AA^2=\spec \CC[u,v]$ in the origin, $E$ the exceptional divisor, $D_1$ the strict transform of $\{u=0\}$, and $D_2$ the strict transform of $\{u+v^2+v=0\}$. Then $$\D=\{\sfrac{2}\}\otimes D_1+\{\sfrac{3}\}\otimes D_2 + [0,\sfrac{6}]\otimes E$$ is a p-divisor with $X=X(\D)$.

Now, although $Y$ is in fact toric, the p-divisor $\D$ isn't invariant with respect to any $\CC^*$-action. Furthermore, the pair $(Y,\D)$ is minimal in the sense that $\D$ is not the pullback of any other p-divisor along some blowup. It follows that $(Y,\D)$ is in fact the pair from the previous remark, so any residual torus action on $X$ would also have to act on $(Y,\D)$. Thus, we can conclude that $X$ doesn't admit any complexity-one action.
\end{ex}

\section{Invariant Semiample Divisors}\label{sec:ssa}
In \cite{petersen:08a}, invariant Cartier divisors on a complexity-one $T$-variety are described via piecewise-affine functions on the slices of a corresponding divisorial fan. Many properties of divisors can be read off directly from these functions. In higher complexity, the situation becomes more convoluted. However, we can still say something if we make some simplifying assumptions.  We first restrict attention to $T$-varieties coming from contraction-free divisorial fans, and then further restrict attention to so-called semicomplete $T$-varieties.

\subsection{Contraction-Free Divisorial Fans}
For the rest of the subsection, we will consider some contraction-free divisorial fan $\dfan$ on
$Y$. We are interested in divisors on the $T$-variety $X=X(\dfan)$. To begin, we associate some additional data with an invariant $\QQ$-Cartier divisor
$$
D=\sum_\rho a_\rho D_\rho+\sum_{P,v}\mu(v)b_{P,v}D_{P,v}.
$$

\begin{lemma}\label{lemma:supportfn}
The $T$-invariant divisor $D \in \CDiv_\QQ(X)$ uniquely determines
piecewise-affine functions $h_P^D\colon|\dfan_P|\to\QQ$ for each prime divisor
$P \subset Y$, satisfying the following conditions.
\begin{enumerate}
\item $h_P^D$ is affine on each $\Delta \in \dfan_P$.
\item For any $\rho\in\xray(\dfan)$, $\slope_\rho(h_P) = -a_\rho$.
\item For any $v\in\xvert_P(\dfan)$, $h_P^D(v)= -b_{P,v}$.
\end{enumerate}
Here, $\slope_\rho(h_P)$ is the slope of $h_P$ along any
one-dimensional polyhedron in $\dfan_P$ with tailcone $\rho$.
\end{lemma}

\begin{proof}
Note that $h_P^D$ is determined by its values at the vertices and along the rays
of $|\dfan_P|$. Since $\dfan$ is contraction-free, these values are prescribed,
so uniqueness is immediate, and it remains to show that affine functions with
these values exist for each $\Delta \in \dfan_P$.

We can assume that $D$ is Cartier, since if the statement is true for $l\cdot 
D$, it will be true for $D$. Since $D$ is $T$-invariant, it is defined by some 
open affine invariant cover $\{U_i\}$ of $X$ together with invariant functions 
$g_i=f_i\cdot\chi^{u_i}$, where $f_i\in\CC(Y)$, see Lemma \ref{lemma:invcover}. 
We consider the divisorial fan induced by the open invariant subsets 
$X(\D)\cap U_i$ for $\D\in \dfan$; this has the same slices as $\dfan$.
Now any $\Delta\in\dfan_P$ which is maximal with respect to inclusion appears as the $P$-coefficient for the p-divisor of some $X(\D)\cap U_i$.
We define an affine function $h \colon \Delta \to \QQ$ by
$h(v) = -\pair{v}{u} - \ord_P f_i$.
Lemma~\ref{lemma:hauptdivisor} shows that for each vertex $v$ of $\Delta$
and for each ray $\rho$ in the tailcone of $\Delta$, we have
\begin{align*}
  b_{P,v} &= \pair{v}{u} + \ord_P f_i = -h(v) \\
  a_\rho  &= \pair{v_\rho}{u} = -(h(v+v_\rho) - h(v)),
\end{align*}
completing the proof.
\end{proof}

These functions are related to the support functions occurring in~\cite{petersen:08a}.
Note that $h_P^D$ determines the function $\Psi_P^D$ introduced earlier. The 
converse is not true in general, but the following characterization 
shows that it does hold for base-point free divisors.
For any point $y\in Y$, let $\mcP(y)$ be the set of all prime divisors passing through $y$. 

\begin{theorem}\label{thm:gg}
Consider a $T$-invariant Cartier divisor $D$ on the $T$-variety $X(\dfan)$ with $\dfan$ contraction-free. Then $D$ is base-point free if and only if for every $\D\in \dfan$ and $y\in \loc \D$
 there exist $u\in \Box^D\cap M$ and $s\in L(\Psi^D(u))$ satisfying
\begin{enumerate}
\item $\Psi_P^D(u)+\ord_P s=0$ for all $P\in\mcP(y)$;
\item $h_P^D(v) = \pair{v}{u} - \Psi_P^D(u)$ for all $P\in\mcP(y)$ and $v\in\D_P$.
\end{enumerate}
\end{theorem}

\begin{proof}
Note that by Remark~\ref{rem:invgen}, a point $x \in X(\dfan)$ is a base-point 
of $D$ if and only if for every homogeneous section $s\cdot\chi^u \in L(D)$, 
$x$ lies in the support of the associated effective invariant divisor
$\Div(s \cdot \chi^u) + D$.

Suppose $D$ is base-point free, and consider any $\D \in \dfan$ and
$y \in \loc \D$. The torus $T$ acts on the fiber $X_y \subset X(\D)$,
and the orbits of this action correspond to the faces of
$\D_y = \sum_{P \in \mcP(y)} \D_P$, see~\cite{altmann:06a}, Section~7.
We choose $x \in X_y$ in a closed orbit of this action, corresponding
to the maximal face of $\D_y$. This point is contained in every vertical
Weil divisor $D_{P,v}$ for $P \in \mcP(y)$, as well as every
horizontal divisor $D_\rho$ for any ray $\rho$ in the tailcone of $\D$.

Now $x$ is not a base-point of $D$, so we find $s \cdot \chi^u \in L(D)$
such that $x \not\in \supp \Div(s \cdot\chi^u) + D$, hence the coefficients
of the invariant divisors through $x$ must cancel out.
By Lemma~\ref{lemma:hauptdivisor}, this means that
\begin{align*}
\pair{v_\rho}{u} + a_\rho &= 0 & \pair{v}{u} + \ord_P(s) + b_{P,v} &= 0.
\end{align*}
for all rays and vertices in the polyhedra $\D_P$.
Applying Lemma~\ref{lemma:supportfn}, it follows that
$$
h_P^D(v) = \pair{v}{u} + \ord_P(s)
$$
on $\D_P$ for all $P \in \mcP(y)$.

By Proposition~\ref{prop:gs}, we know that $s$ is a section of $\Psi^D(u)$,
so $\Psi_P^D(u)+\ord_P(s)\geq 0$ for all $P$. On the other hand,
by definition of $\Psi_P^D$, we have
$$
\Psi_P^D(u) \leq \pair{v}{u} - h_P^D(v) = -\ord_P(s)
$$
for all $P\in\mcP(y)$ and $v\in\D_P$. 
It follows that $\Psi_P^D(u) = -\ord_P(s)$, so 
points one and two of the theorem are satisfied.

Conversely, let us assume that $D$ satisfies the conditions of the theorem, and 
show that $D$ is base-point free. Take any $x \in X$ mapping to $y \in Y$, lying 
in $X(\D)$ for some $\D \in \dfan$. With $u$ and $s$ as given by the hypothesis, 
it suffices to show that $\Div(s\cdot\chi^u) + D$ doesn't meet $X_y$. But this 
follows by inverting the application of Lemmas~\ref{lemma:hauptdivisor}
and~\ref{lemma:supportfn} in the first part of the proof,
since we know that $h_P^D(v) = \pair{v}{u} + \ord_P(s)$.
\end{proof}

We now draw some consequences from the above theorem. Consider any concave piecewise-affine function $\Psi\colon\Box\to\WDiv_{\QQ\cup\{\infty\}} Y$, that is the minimum of only finitely many affine functions.
The \emph{lineality space} of $\Psi_P$ is the largest subspace $L_P\subset M_\QQ$ such that $\Psi_P$ and $\Box$ are invariant under translation by elements of $L_P$. Thus, we can considered $\Psi_P$ to be defined on $\Box/L_P$.
The \emph{vertices} of $\Psi_P$ are those $\overline u\in\Box/L_P$ such that $(\overline u,\Psi_P(\overline u))$ is a vertex of the graph of $\Psi_P$.
\begin{definition}\label{def:sharp}
We say that $\Psi$ is \emph{asymptotically sharp} if for any prime divisor $P\subset Y$ with $\Psi_P\not\equiv\infty$ ,  the following holds for all vertices $\overline u$ of $\Psi_P$: 
\begin{equation}\tag{$**$}\label{eqn:sharp}
\begin{array}{l}
\textrm{There exists $k>0$ and $u\in \Box\cap (\overline u+L_P)$ such that there is some $s\in L(k \Psi(u))$}\\
\textrm{satisfying $\ord_Ps+k \Psi_P(u)=0$. 
}
\end{array}
\end{equation}
\noindent{}We say that $\Psi$ is \emph{sharp} if we can always take $k=1$ and $u\in M$.
\end{definition}

\begin{cor}\label{cor:sa}
Let $D$ be a $\QQ$-Cartier divisor on a  $T$-variety $X(\dfan)$, where $\dfan$ is contraction-free and $|\dfan_P|$ is convex  for all prime divisors $P$.  Then if $D$ is base-point free/semiample, it follows that 
\begin{enumerate}
\item $\Psi^D$ is sharp/asymptotically sharp; and
\item For any $P\subset Y$ prime, $h_P^D$ is concave.
\end{enumerate}
\end{cor}
\begin{proof}
We will prove the statement for the base-point free case.  The semiample case follows immediately by passing to a sufficient multiple of $D$.
We first remark that $h_P^D$ is concave if and only if for all $\Delta \in \dfan_P$, there is $u\in\Box^D$ such that $\Psi^D(u)=\pair{v}{u} -h_P^D(v)$ for any $v\in \Delta$.

Now, consider any prime divisor $P\subset Y$. Then by the remark above and Theorem \ref{thm:gg}, $h_P^D$ must be concave. For any maximal dimensional $\Delta\in\dfan_P$, $h_P^D(v)$ restricted to $\Delta$ is of the form $\pair{v}{\overline u} + \ord_P s$ for some uniquely determined $\overline u\in M/L_P$, which must be a vertex of $\Psi_P^D$. Furthermore, all vertices of $\Psi_P^D$ arise in this way. This together with the above theorem implies the sharpness condition.  
\end{proof}

We now draw another consequence of our characterization of global generation.
\begin{definition}\label{def:divpolysum}
Consider any concave piecewise-affine maps $\Psi^i\colon\Box^i\to\WDiv_{\QQ\cup\{\infty\}} Y$, $i=1,2$. We define their sum to be
\begin{align*}
(\Psi^1+\Psi^2)\colon(\Box^1+\Box^2)&\to\WDiv_{\QQ\cup\{\infty\}} Y\\
u&\mapsto \sum_P \max_{\substack{u_i\in\Box^i\\u_1+u_2=u}} (\Psi_P^1(u_1)+\Psi_P^2(u_2))P.
\end{align*}
\end{definition}
\begin{prop}
Let $D$ and $E$ be semiample $T$-invariant $\QQ$-Cartier divisors on
$X$. Then $\Box^D+\Box^E=\Box^{D+E}$ and
$
\Psi^D+\Psi^E=\Psi^{D+E}
$.
\end{prop}
\begin{proof}
This is a consequence of Corollary \ref{cor:sa}. Indeed, since
$h_P^D,h_P^E$ are concave, the claim follows from 
general facts about convexity.
\end{proof}

\subsection{Semicompleteness}
We will call a variety $X$ \emph{semicomplete} if the
affine contraction morphism $r \colon X \to X_0 = \Spec H^0(X, \CO_X)$ 
is proper. Equivalently, a variety $X$ is semicomplete if it is proper over any affine variety. This is closely related to semiprojectivity:

\begin{prop}\label{prop:semicomplete}
A variety $X$ is semiprojective if and only if it is quasiprojective and semicomplete.
\end{prop}
\begin{proof}
If $X$ is semiprojective, then the claim is immediate. Instead, assume that $X$ is quasiprojective and semicomplete.
Thus $X$ is embedded in some $\PP^m$ and we have a proper map $\theta\colon X\to Z$ for some affine $Z$. Then the image of the graph $\Gamma$ of $\theta$ in $\PP^m\times Z$ is isomorphic to $X$, and furthermore, we claim that it is closed. Indeed, we have maps $\Gamma \to \PP^m\times Z \to Z$, with the first separated, and the second and the composition proper. Thus, the first is proper as well and the image of $\Gamma$ in $\PP^m\times Z$ is closed, making $X$ semiprojective.
\end{proof}

As before, let $X = X(\dfan)$ be a $T$-variety, with $\dfan$ a divisorial fan on $Y$.
The affine contraction $X_0$ retains an action by $T$ which need not
be effective. There is an effective residual action by a 
quotient torus $T'$, corresponding to the sublattice $M'$ of
$M$ generated by the weight cone of the $T$-action.
In order to be able to analyze the equivariant morphism $X \to X_0$,
we recall some facts and definitions about morphisms of polyhedral
divisors.

\begin{definition}\label{def:map} Let $\D$ be a $\QQ$-Cartier polyhedral divisor on some $Y'$ with coefficients in $N'$  and $\dfan$ a divisorial fan on $Y$ with coefficients in $N$.
\begin{enumerate}
\item The pull-back of $\D$ under a dominant morphism
	$\psi \colon Y \to Y'$ is the map $\psi^*\D\colon\wt\to \WDiv_{\QQ\cup\{\infty\}}(Y)$ given by
      $$
      (\psi^*\D)(u) = \psi^*(\D(u)). 
      $$
\item The pull-back of $\D$ under $F \colon N \to N'$ is defined
      by
      $$
      F^{-1}(\D)_P = F^{-1}(\D_P).
      $$
\item The principal polyhedral divisor associated with a function
      $\mathfrak{f}=\sum_i v_i\otimes f_i \in N \otimes \CC(Y)^*$ is defined by
      $$
      \Div(\mathfrak{f})(u) = \sum_i\pair{v_i}{u} \Div(f_i).
      $$
\item The pull-back of $\D$ under a triple
      $\varphi = (\psi, F, \mathfrak{f})$ is
      $$
      \varphi^{-1}(\D) = F^{-1}(\psi^*(\D) + \Div(\mathfrak{f})).
      $$
\item A \emph{map} $\dfan \to \D$ is a triple $\varphi$ as above,
      such that $|\dfan_P| \subset \varphi^{-1}(\D)_P$ for all prime
      divisors $P$ on $Y$.
\end{enumerate}
\end{definition}

\begin{theorem}[cf. {\cite[Section 8]{altmann:06a}}]\label{thm:morphisms} 
A map $\varphi \colon \dfan \to \D$ defines an equivariant morphism
$\varphi \colon X(\dfan) \to X(\D)$.
If $\varphi \colon X(\dfan) \to X(\D)$ is an equivariant morphism and $\D$ is proper,
there exists a modification $p \colon \widetilde{Y} \to Y$ and a
map $\varphi' \colon p^*\dfan \to \D$ defining
$\varphi \colon X(\dfan) = X(p^*\dfan) \to X(\D).$
\end{theorem}

We can now provide a polyhedral description of the equivariant morphism
$X \to X_0$.

\begin{lemma}\label{lemma:psi0}
Assume that $\dfan$ is contraction-free with tailfan $\Sigma$.
\begin{enumerate}
\item
  The function $\Psi^0 \colon \Box^0 \to \WDiv_{\QQ\cup\{\infty\}}(Y)$ associated
  with the trivial divisor on $X$ is a polyhedral divisor
  with tailcone $(\Box^0)^\vee = \sigma := \conv |\Sigma|$ and 
  coefficients $\Psi^0_P = \conv |\dfan_P|$. 
\item
  The algebra of global sections is the coordinate ring of $X_0$,
  that is, $X_0 = X(\Psi^0)$.
\item
  $r \colon X \to X_0$ is defined by
  the map $(\id_Y, \id_N, 1) \colon \dfan \to \Psi^0$.
\end{enumerate}
\end{lemma}

\begin{proof}
For the first claim, note that by definition,
$\Box^0$ is a polyhedral cone and $\Psi^0$ is convex and
piecewise-linear. The rest follows by dualizing the definitions of
$\Box^0$ and $\Psi^0$ and using that  $\dfan$ is contraction-free. 

The second claim is a direct consequence of Proposition~\ref{prop:gs}.
The third follows from the definition of the morphism associated with
a map of polyhedral divisors.
\end{proof}

In general, $\Psi^0$ is not a p-divisor for $X_0$. But if $X$ is
semicomplete, we will see that it is semiample. We need the
following result on properness of maps of polyhedral divisors.

\begin{theorem}[cf. {\cite[Satz 2.16]{suess:10a}}]\label{thm:proper}
If a map $\varphi \colon \dfan \to \D$ defines a proper morphism for $\dfan$ a divisorial fan and $\D$ a p-divisor, then for all prime divisors $P\subset Y$, $|\dfan_P| = \varphi^{-1}(\D)_P$.
\end{theorem}

\begin{prop}\label{prop:semiproj}
If $X=X(\dfan)$ is semicomplete for a contraction-free $\dfan$ and $\Psi^0$ is $\QQ$-Cartier, then $\Psi^0$ is semiample, and for all prime divisors $P$ on $Y$,
$|\dfan_P|$ is convex. Furthermore, $\loc \Psi^0$ is semicomplete.
\end{prop}

\begin{proof}
Let the p-divisor $\D$ on $Y'$ describe the $T'$-variety $X_0$.
By Proposition~\ref{thm:morphisms}, the equivariant morphism
$X \to X_0$ corresponds to a map $\dfan \to \D$,
blowing up $Y$ if necessary. The pulled back $\dfan$ is not necessarily contraction-free anymore, so we replace it by $\dfan'$ as in Proposition~\ref{prop:contraction-free}.
Then the affine contractions of $X(\dfan)$ and $X(\dfan')$ agree, and the map $r \colon X(\dfan')\to X_0$ corresponds to a map $\varphi \colon \dfan' \to \D$.

Since $X$ is semicomplete, $r$ is proper, so by
Theorem~\ref{thm:proper}, we get $|\dfan_P|=|\dfan'_P|=\varphi^{-1}(\D')_P$. 
This shows the convexity claim. Now by Lemma \ref{lemma:psi0}, we have that  $|\dfan_P|=\Psi_P^0$. Thus, the pullback of $\Psi^0$  is the pullback of a p-divisor (which is in particular semiample), so it must be semiample itself.

To show the semicompleteness of $\loc \Psi^0$, we first note that the map $\loc \dfan' \to \loc \D$ is proper; this follows from \cite[Satz 2.16]{suess:10a}. From the semiprojectivity of $\loc \D$, we thus have a proper map $\theta'\colon \loc \dfan'\to Z$, where $Z$ is affine. Furthermore, $\theta'$ factors through $\loc \dfan$, since the regular functions on $\loc \dfan$ and $\loc \dfan'$ are equal, and $Z$ is affine. Let $\theta$ denote this map from $\loc \dfan$ to $Z$. We claim that $\theta$ is proper. Indeed, this follows from the separatedness of $\theta$, the surjectivity of $\loc \dfan'\to\loc\dfan$, and the properness of $\theta'$, see \cite[Corollary 5.4.3]{EGA2}. 
\end{proof}

\begin{ex}\label{ex:nsa}
Let $Y$ be the blowup of $\AA^2$ at the origin; let $D_1$, $D_2$
be the strict transforms of the coordinate axes and $E$ the exceptional divisor. Consider the
divisorial fan $\dfan$ generated by
{$$
\setlength{\arraycolsep}{2pt}
\begin{array}[b]{rcrcrcr}
\D^1 &=& [0,\infty) \otimes D_1 &+& [1,\infty) \otimes E &+& \emptyset \otimes D_2 \\
\D^2 &=& \emptyset  \otimes D_1 &+& [1,\infty) \otimes E &+& [0,\infty) \otimes D_2
\end{array}.
$$}
Then $X$ is not semicomplete, and $\Psi^0$ is big but not semiample:
$$
\Psi^0 = [0,\infty) \otimes D_1 + [1,\infty) \otimes E + [0,\infty) \otimes D_2,
$$
so $\Psi^0(1) = E$. A p-divisor for $X_0 = \AA^3$ is
$$
\D^0 = [0,\infty) \otimes D_1 + [0,\infty) \otimes D_2,
$$
which may be defined on $Y_0 = \Spec H^0(Y, \CO_Y) = \AA^2$.
\end{ex}

\begin{ex}
Consider the same example as above, but modify $D_1$ and $D_2$ so that
$E$ has coefficient $[-1, \infty)$. Then $\Psi^0$ is a p-divisor, and
$X$ is semicomplete. In contrast to the previous example, $X_0$ can
not be expressed by a p-divisor on $Y_0$.
\end{ex}

\begin{ex}
Let $X = \PP^1 \times \PP^1$, with a complexity one action through
one factor. The divisorial fan for $X$ is defined on $Y = \PP^1$.
Then $\Psi^0$ is a semiample divisor on $\PP^1$ for the point $X_0$,
but it is not big.
\end{ex}

\subsection{Semiample Decomposition}
Let $X(\dfan)$ be some $T$-variety, and $D$ an invariant divisor on $X(\dfan)$.
If $D$ is semiample, it  does not in general hold that $\Psi^D(u)$ is also semiample for all $u\in\Box^D$. We wish to address this question in the situation where $\dfan$ is contraction-free and $X(\dfan)$ is semicomplete. We first need the following lemma:

\begin{lemma}\label{lemma:sadecomp}
Let $\dfan$ be contraction-free with $\dfan_P$ convex and $\loc \dfan$ semiprojective. We consider the variety $X=X(\dfan)$ and assume that  $\Psi^0(u)$ is semiample for all $u\in\Box^0$. Consider any semiample invariant divisor $D$ on $X$. 
If there is a concave semiample-valued $\Psi'\colon\Box^D\to\WDiv_{\QQ\cup \{\infty\}} Y$ with semiprojective locus such that  $L(k\Psi^D(u))=L(k \Psi'(u))$ for all $u\in \Box^D\cap \sfrac{k}M$ and all $k\in\ZZ_{\geq 0}$,  then $\Psi^D(u)$ is semiample for all $u\in\Box^D$.
\end{lemma}

\begin{proof}
Let $\Psi'$ be as in the hypothesis. It follows from Lemma 9.1 of \cite{altmann:08a} that $\Psi_P^D\geq \Psi_P'$ as long as $\Psi_P'\not\equiv \infty$. Note that the proof of the lemma does not require that the second divisor is semiample.
Furthermore, if $\Psi_P'\equiv \infty$, then we must have $\Psi_P^D\equiv \infty$ as well. Indeed, the complement of $\loc \Psi'$ is some semiample divisor $C$. Thus, any nonempty $L(k\Psi'(u))$ has a section with an arbitrarily large pole along $C$. Similarly, if $\Psi_P^D\equiv \infty$, then $\Psi_P'\equiv \infty$ due to the semiprojectivity of $\loc \dfan$.

Fix now some prime $P$ with $\dfan_P\neq \emptyset$.
Let $\overline{\Box}_P\subset \Box^D/L_P$ be the convex hull of the vertices of $\Psi_P$ and $\Box_P$ its inverse image in $\Box^D$. It follows from asymptotic sharpness that for any vertex $\overline u$ of some $\Psi^D_P$, there is some $u\in\Box^P$ mapping to $\overline u$ with $\Psi_P^D(u)=\Psi_P'(u)$.
Consider any $w\in L_P$.
Then we even have $\Psi_P^D(u+w)=\Psi_P'(u+w)$. Indeed, $\Psi^D(u+w)\geq \Psi^D(u)+\Psi^0(w)$, and some multiple $k$ of the left hand side thus has a section vanishing along $P$ of order $$k(\Psi_P^D(u)+\Psi_P^0(w))=k\Psi_P^D(u+w)$$ by choice of $u$ and the semiampleness of $\Psi^0(w)$. It even follows by convexity that $\Psi_P^D(u)=\Psi_P'(u)$ for any $u\in\Box_P$.

Now consider any $u\in\Box_P$ and $w\in\tail(\Box^D)$. Then there exists $l \gg 0$ such that for $\lambda\geq 0$
\begin{align*}
\Psi^D(u+(l+\lambda)w)=\Psi^D(u+lw)+\Psi^0(\lambda w);\\
\Psi'(u +(l+\lambda) w)\leq\Psi'(u+lw)+\Psi^0(\lambda\cdot w).
\end{align*}
Since the right hand side of the second line above is semiample, we must actually have equality again by Lemma 9.1 of \cite{altmann:08a}.
The concavity of $\Psi'$ and equality of $\Psi_P^D$ and $\Psi_P'$ on $\overline\Box_P$ together with the above imply
$$
2\Psi_P'(u+lw)\geq \Psi_P'(u+2lw)+\Psi_P'(u)=
\Psi_P'(u+lw)+\Psi_P^0(lw)+\Psi_P^D(u).
$$
We can thus conclude that $\Psi_P'(u+lw)\geq \Psi_P^0(lw)+\Psi_P^D(u)=\Psi_P^D(u+lw)$. From the concavity of $\Psi'$ we then get $\Psi_P^D(u+\lambda w)=\Psi_P'(u+\lambda w)$ for any $\lambda\geq 0$.  But any $u'\in\Box^D$ can be written as such a sum $u+\lambda w$. Thus, $\Psi^D=\Psi'$, so $\Psi^D$ is semiample-valued.
\end{proof}

\begin{theorem}\label{thm:sadecomp}
Let $\dfan$ be a contraction-free divisorial fan on smooth $Y$ with $X=X(\dfan)$ semicomplete. Consider any semiample divisor $D$ on $X$. Then $\Psi^D(u)$ is semiample for all $u\in\Box^D$. Furthermore, if $D$ is big, then $\Psi^D(u)$ is big for all $u\in \relint \Box^D$. 
\end{theorem}
\begin{figure}
$$
\xymatrix{
 X(\rho^*\D) \ar@{-->}[d] \ar@{=}[r] &X(\D)\ar@{-->}[d]\ar@{=}[r] &X(\tD)\ar@{-->}[dd]\\
 X(\dfan')\ar[dd] \ar[r]^\rho&X(\dfan) \ar[d] &\\
 & Y & Z\\
W \ar[ur]\ar[urr]        }
$$
\caption{Situation in proof of Theorem \ref{thm:sadecomp}}\label{fig:sadecomp}
\end{figure}

\begin{proof}
Without loss of generality, we can assume that $Y$ is complete. Indeed, we can complete $Y$ and pull back the divisorial fan without changing it from being contraction-free.

Now, let $\D$ be a $T$-invariant p-divisor on $X$ with some weight cone $\wt\subset \ZZ^2_\QQ$, such that for some $w_0\in\wt$, $\D(w_0)=D$; such a $\D$ always exists. If $D$ is big, we can even require $w_0\in\relint\wt$. Then $X(\D)$ is a $(\CC^*)^2$-variety. But it also inherits the $T$-action of $X(\dfan)$, so it is in fact a $\tT$ variety, where $\tT=\CC^*\otimes(N\oplus\ZZ^2)$.

Let $Z$ be a projective completion for the normalization of the special component of the Chow quotient of $X(\D)$ by the action of $\tT$, see Section 6 of \cite{altmann:06a}. Let $\tD$ be the corresponding p-divisor on $Z$ with $X(\tD)=X(\D)$. One easily checks that $Z$ and $Y$ are birational; we thus find a projective variety $W$ mapping properly and birationally to both $Y$ and $Z$.
We now pull back $\dfan$ to $W$, possibly blowing up to make it contraction-free again, giving us a new divisorial fan $\dfan'$ and a map $\rho\colon X(\dfan')\to X(\dfan)$. The situation thus far is pictured in Figure \ref{fig:sadecomp}.

Now, since $X$ is semicomplete, $\Psi^0$ is semiample, $\loc \dfan$ is semicomplete, and $\dfan_P$ is convex by Proposition~\ref{prop:semiproj}. Furthermore, since $\loc(\dfan)$ was semicomplete, $\loc(\dfan')$ will be semicomplete as well, and thus semiprojective by Proposition~\ref{prop:semicomplete}.

Note that for any invariant Cartier divisor $E$ on $X(\dfan)$, $\Psi^{\rho^* E}$ is simply the  pullback of $\Psi^E$ to $W$.  We now pull back $\tD$ to $W$, and after possibly correcting with some principal polyhedral divisor, we have $L(\Psi^{\D(w_0)}(u))=L(\tD(u,w_0))$ for all $u\in\Box^{\D(w_0)}\cap M$ by Proposition~\ref{prop:gs}. Setting $\Psi'(u)=\tD(u,w_0)$, we can apply the above lemma and conclude that the pullback of $\Psi^D(u)$ to $W$ must have been semiample. Thus, $\Psi^D(u)$ must have been semiample as well. Furthermore, from the proof of the above lemma, we actually have that the pullback of $\Psi^{\D(w_0)}(u)$ is the pullback of $\Psi'(u)$, which is big for $u\in\relint \Box^D$ if $w_0\in\relint \wt$. Thus, if $D$ is big, $\Psi^{\D(w_0)}(u)$ is big as well.
\end{proof}

\begin{rem}
As we can see by Example \ref{ex:nsa}, the above theorem does not necessarily  hold if $X(\dfan)$ is not semicomplete.
\end{rem}

\subsection{Upgrades revisited}\label{sec:upgradesrev}
We are now in a position to prove Theorem \ref{thm:tD}:
\begin{proof}[Proof of Theorem \ref{thm:tD}.]
It follows directly from the dual description of 
Proposition~\ref{prop:dualupgrade} that $\twt$ is a polyhedral cone
and that $\tD$ is a polyhedral divisor. The fact that $X=X(\tD)$ follows easily from Proposition~\ref{prop:gs}. If $Y$ is semiprojective and $\dfan$ is contraction-free, then $\loc \tD$ is semiprojective by Propositions~\ref{prop:semicomplete} and~\ref{prop:semiproj}. Furthermore, we can apply Theorem \ref{thm:sadecomp} to prove that $\tD$ is semiample, and big on the interior of $\twt$.
The cone $\twt$ must be full-dimensional since $\dim X=\dim \twt + \dim Y'$.
\end{proof}

\section{Divisorial Polyhedra}\label{sec:divpoly}
In Section \ref{sec:downgradeone}, we will be considering the problem of downgrading the torus action of a complexity-one $T$-variety. To prepare, we need a method for construction semiprojective complexity-one $T$-varieties. We do this by adapting the so-called divisorial polytopes introduced by H. S\"u\ss{} and the first author in \cite{ilten:09d}.
Throughout this section, $Y$ will be a smooth curve. As usual, $M$ is some lattice with dual $N$.
\begin{definition}
	A \emph{divisorial polyhedron} consists of a pair $(\Psi,\Box)$, where $\Box$ is a polyhedron in $M_\QQ$, and $\Psi$ is a piecewise-affine concave map from $\Box$ to $\CDiv_{\QQ} Y$ taking values in semiample divisors. Often we will refer to such a divisorial polyhedron simply by $\Psi$. Note that we can add two divisorial polyhedra just as in Definition \ref{def:divpolysum}.
\end{definition}

Given a divisorial polyhedron $\Psi$, we will show how to associate a semiprojective complexity-one $T$-variety and a $\QQ$-Cartier $\QQ$-divisor.
For $P\subset Y$ prime, we define its \emph{linear part}  as follows:
\begin{align*}
\Psi_P^\lin\colon\tail(\Box)&\to \QQ \\
w&\mapsto \lim_{\lambda\to\infty} \Psi_P(u+\lambda w)/\lambda
\end{align*}
for any $u\in\Box$. We then associate dual polyhedra $\Box_P^*$ in the $N$ lattice as well as functions $\Psi_P^*\colon\Box_P^*\to\QQ$:
\begin{align*}
	\Box_P^*&:=\{v\in N_\QQ\ |\ \pair{v}{w} \geq \Psi_P^\lin(w)\ \forall\ w\in\tail(\Box)\};\\
	\Psi_P^*(v)&:=\min_{u\in \Box} (\pair{v}{u} -\Psi_P(u)).
\end{align*}
We let $\Psi^*$ be the formal sum $\sum \Psi_P^*\otimes P$.
For any piecewise-affine function $f$, denote by $\Xi(f)$ the polyhedral subdivision induced on its domain. Now assume that $\Box$ is full-dimensional. Then $\Xi(\Psi_P^*)$ is a subdivision of $\Box_P^*$ consisting of pointed polyhedra for any $P$. 

We are now ready to construct a divisorial fan. Let $\mcP$ be some finite nonempty set of points $P$ in $Y$ containing all $P$ with $\Psi_P\not\equiv 0$ and consider the divisor $E:=\sum_{P\in\mcP} P$.
Let  $\C$ be the set consisting of p-divisors on $Y$ of the form
$$\Delta_P\otimes P + \emptyset \otimes (E-P)$$
for any $P\in\mcP$  and any $\Delta_P\in\Xi(\Psi_P^*)$.

\begin{prop}\label{prop:divpoly}
$\C$ induces a contraction-free divisorial fan $\dfan$ via intersection. The corresponding variety $X:=X(\dfan)$ is semiprojective, and the support function $\Psi^*$ corresponds to a semiample $\QQ$-divisor on $X$.
\end{prop}
\begin{proof}
To see that $\C$ induces a contraction-free divisorial fan, we just need to check that $\tail(\Xi(\Psi_P^*))$ is independent of $P$. Likewise, to check that $\Psi^*$ corresponds to a divisor on $X$, we need to show the independence from $P$ of 
\begin{align*}
(\Psi_P^*)^\lin\colon\tail(\Box_P^*)&\to \QQ \\
w&\mapsto \lim_{\lambda\to\infty} \Psi_P(v+\lambda w)/\lambda
\end{align*}
for any $v\in\Box_P^*$, see \cite[Section 3.1]{petersen:08a}. This would of course imply that $\tail(\Xi_P^*)$ is independent of $P$, so we shall just check the latter. But a straightforward calculation shows that $(\Psi_P^*)^\lin$ is just the map $\tail(\Box)^\vee\to\QQ$ sending $w$ to $\min_{u\in\Box}\pair{w}{u}$.

We now show that we can embed $X$ as an open subvariety of a projective variety $\overline X$. Let $\overline Y$ be the unique smooth projective curve birational to $Y$ (already equal to $Y$, if $Y$ was projective). Let $\Box'$ be a full-dimensional polytope such that $\Box=\Box'+\tail\Box$ and containing at least one point from the relative interior of each domain of linearity of $\Psi$. We take $\Psi'$ to be the divisorial polyhedron $\Box'\to \CDiv_\QQ \overline Y$ with $\Psi_P'$ the restriction of $\Psi_P$ for $P\in Y$, and $\Psi_P'=0$ for $P\in\overline{Y}\setminus Y$.
Note that from $\Box$ and $\Psi'$, we can reconstruct $\Psi$.
In any case, from $\Psi'$, we get a divisorial fan $\dfan'$; $\overline X=X(\dfan')$ is projective by a slightly modified version of \cite[Theorem 3.2]{ilten:09d}. Likewise, the divisor corresponding to $(\Psi')^*$ is semiample. 

We claim that there is an open embedding of $X$ in $\overline X$. Indeed, consider any $\Delta\in\Xi(\Psi_P^*)$. This is also an element of $\Xi( (\Psi')_P^*)$. Thus, $\dfan_P$ is a subcomplex of $\dfan_P'$ for every $P$ and the claim follows readily. Furthermore, $\Psi^*$ is the restriction of $(\Psi')^*$ to $\dfan$, so the semiampleness of the divisor corresponding to the support function $\Psi^*$ follows.

It remains to be shown that $X$ is semiprojective. However, it will be enough to show that $X$ is proper over something affine since we have seen that it is quasiprojective, see Proposition~\ref{prop:semicomplete}. Consider the map $r\colon X\to\X_0$ from $X$ to its affine contraction $X_0$. Let $\D=\sum \Box_P^*\otimes P$ be a polyhedral divisor on $Y$. Note that $X_0=X(\D)$.
Let $L$ be the largest linear subspace contained in $\tail(\D)$ and consider $F\colon N\to N':=N/L\cap N$. Set $\D':=\sum F(\D_P)\otimes P$.
If $Y\neq\overline Y$, then $\D'$ is a p-divisor.
The properness of $r$ then follows immediately from \cite[Satz 2.16]{suess:10a}, applied to the map $\dfan\to \D'$ coming from $(\id_Y,F,1)$.

We thus assume that $Y=\overline Y$. It follows that $\deg \D'(u)\geq 0$ for all $u\in\tail(\D')^\vee$. If $\D'$ is big, then it is in fact a p-divisor, since $X(\D')$ has a fiitely generated coordinate ring. Similar to above, $r$ must be proper.
If $\D'$ isn't big, $\deg \D'(u)=0$ for all $u$. We must thus have $\D'=\sum (v_i+\tail \D')\otimes \Div f_i$ for some $v_i\in N'_\QQ$, $f_i\in\CC(Y)$. Let $G$ be the inclusion $N'\to N'':=N[v_i]$ and $M''$ the corresponding sublattice of $M'$; all the weights $u$ such that $\D'(u)$ has a section actually lie in $M''$. Let $\mathfrak{f}=\sum G(\D_P')\otimes P$; this is a principle polyhedral divisor for $N''$. The properness of $r$ then follows again from \cite[Satz 2.16]{suess:10a}, this time considering the map $\dfan\to \tail(\mathfrak{f})$ gotten by the composition $(Y\to \pt,\id_{N''},\mathfrak{f})\circ (\id_Y,G\circ F,1)$.
\end{proof}

\section{Downgrading a Torus Action}\label{sec:downgrade}
We now briefly turn our back on the problem of upgrading a torus action, and instead consider the problem of downgrading a torus action, i.e. restricting to some subtorus. We shall first solve this problem for complexity-one $T$-varieties, and then discuss problems arising for downgrades of higher complexity.

\subsection{Complexity One}\label{sec:downgradeone}
Let $Y$ be a smooth curve, and $\D$ a p-divisor with locus $Y$ and with coefficients in $N_\QQ$. Let $\pr\colon M\mapsto \hM$ be a surjection of lattices; this corresponds to a closed embedding of the torus $\hT=\hN\otimes\CC^*\hookrightarrow T=N\otimes\CC^*$. We wish to downgrade the torus action and describe $X(\D)$ as a $\hT$ variety. This will involve finding a semiprojective variety $\hY$ of dimension $\dim Y + \rk \ker \pr$ along with a p-divisor $\hD$ with locus $\hY$ such that $X(\D)=X(\hD)$.

Now, we consider the lattice  $M':=\ker \pr$ and choose a cosection $t\colon M\to M'$, inducing a section $s^*\colon \hM\to M$:
\begin{equation*}
\xymatrix{
0 \ar[r]& M' \ar[r]  &M \ar^{\pr}[r] \ar^t@/^1pc/[l] & \hM \ar[r] \ar^{s^*}@/^1pc/[l]& 0.
}
\end{equation*}
 Let $\wt$ be the weight cone of $\D$ and $\hwt$ its image under $\pr$. For $\hu\in\hwt\cap\hM$, we define 
\begin{align*}
\Box[\hu]&:=t(\pr^{-1}(\hu)\cap\wt)\\
\Psi[\hu]&\colon \Box[\hu]\to \CDiv_{\QQ} Y\\
&\quad\quad u'\mapsto\D(u'+s^*(\hu)).
\end{align*}
Now, each $\Psi[\hu]$ is a divisorial polyhedron. For each domain of linearity $\wt_i$ of $\D$, take some $\hu_i$ in the relative interior of $\pr \wt_i$. We let $\dfan$ be a divisorial fan corresponding to $\sum \Psi[\hu_i]$ and $\hY=X(\dfan)$. Note that $\Box[\hu_i]$ is always full-dimensional.
Furthermore, any $\Psi_P[\hu]^*$ is piecewise-affine on each $\dfan_P$, since $\dfan_P$ is the coarsest common subdivision of all the $\Xi(\Psi_P^*(\hu_i))$.

We also consider the dual sequence of lattices with section and cosection:
\begin{equation*}
\xymatrix{
0 \ar[r]& \hN \ar[r]  &N \ar^{\pi}[r] \ar^{s}@/^1pc/[l] & N' \ar[r] \ar^{t^*}@/^1pc/[l]& 0
}
\end{equation*}

For $P\in Y$, we take $\Xi_P$ to be the coarsest polyhedral subdivision of $\pi(\D_P)$ containing subdivisions of $\pi(\Delta)$ for any face $\Delta\prec\D_P$.
For any vertex $v\in \Xi_P$, we take
$$\Delta_{P,v}:=  s(\D_P\cap \pi^{-1}(v)).$$
Likewise, for any ray $\rho\in\tail(\Xi_P)$, we take
$$\Delta_\rho:=s(\tail(\D)\cap \pi^{-1}(v_\rho)).$$

\begin{prop}\label{prop:dcoeff}
In this situation, $\Xi_P=\dfan_P$. Furthermore, for any $\hu\in\hwt$,

$\Psi_P[\hu]^*(v)=\min \pair{\Delta_{P,v}}{\hu}$, and
$(\Psi[\hu]^*)^\lin(v_\rho)=\min \pair{\Delta_\rho}{\hu}$. 
\end{prop}
\begin{proof}
All statements will follow from the duality described in Section 8 of \cite{altmann:03a} for toric downgrades.
To begin with, we have 
\begin{align*}
(\Psi[\hu]^*)^\lin(v_\rho)&=\min \pair{v_\rho}{t(\pr^{-1}(\hu)\cap\wt)} \\
\min \pair{\Delta_\rho}{\hu}&=\min \pair{s(\pi^{-1}(v_\rho)\cap\wt^\vee)}{\hu}
\end{align*}
so $(\Psi[\hu]^*)^\lin(v_\rho)=\min \pair{\Delta_\rho}{\hu}$ follows from the proof of \cite[Proposition 8.5]{altmann:03a}.

For the other claims, we will be using the mutually dual sequences
\begin{align*}
\xymatrix{
0 \ar[r]& M'\oplus\ZZ \ar[r]  &M\oplus\ZZ \ar^{\widehat{\pr}}[r] \ar^{\widehat{t}}@/^1pc/[l] & \hM \ar[r] \ar^{\widehat{s}^*}@/^1pc/[l]& 0
}\\
\xymatrix{
0 \ar[r]& \hN \ar[r]  &N\oplus\ZZ \ar^{\widehat{\pi}}[r] \ar^{\widehat{s}}@/^1pc/[l] & N'\oplus \ZZ \ar[r] \ar^{\widehat{t}^*}@/^1pc/[l]& 0
}
\end{align*}
where $\widehat\pr$ is just the projection $M\oplus\ZZ\to M$ composed with $\pr$, and $\widehat{s}^*$ is $s^*$ composed with the map of $M$ onto the first factor of $M\oplus\ZZ$.

Now for any $P\in Y$, we define a cone $\widehat{\wt}_P\subset (M\oplus\ZZ)_\QQ$:
$$
\widehat{\wt}_P=\{(u,w)\ |\ u\in\wt\ \textrm{and}\ w\geq-\D_P(u)\}.
$$
A straightforward calculation shows $\widehat{\wt}_P^\vee=\pos((\D_P,1))$.
Furthermore, we have 
	\begin{align*}
	\Psi_P[\hu]^*(v)&=\min \pair{(v,1)}{\widehat{t}(\widehat{\pr}^{-1}(\hu)\cap\widehat\wt_P)}\\
\min pair{\Delta_{P,v}}{\hu}&=\min\pair{\widehat{s}(\widehat\wt_P^\vee\cap\widehat\pi^{-1}((v,1)))}{\hu}
\end{align*}
so
$\Psi_P[\hu]^*(v)=\min \pair{\Delta_{P,v}}{\hu}$ again follows from the proof of \cite[Proposition 8.5]{altmann:03a}.

Finally, from Section 8 of \cite{altmann:03a} we also have that the image fan $\widehat\pi(\widehat\wt_P^\vee)$ is equal to the normal fan $\N(\Sigma_i \widehat{t}(\widehat\pr^{-1}(\hu_i)\cap\widehat\wt_P))$. Intersection the former with the hyperplane in height one gives $\Xi_P$, and intersecting the latter gives $\dfan_P$.
\end{proof}	

\begin{theorem}\label{thm:downgrade}
Each support function $\Psi[\hu]^*$ represents a semiample Cartier divisor on $\hY$. The corresponding map $\hD\colon \hwt\to\CDiv_{\QQ} \hY$ is a p-divisor on $\hY$ with $X(\hD)=X(\D)$, and can be written as
$$\hD=\sum_\rho \Delta_\rho\otimes D_\rho+\sum_{P,v} \Delta_{P,v}\otimes \mu(v)D_{P,v}.$$
\end{theorem}

\begin{proof}
The fact that $\Psi[\hu]^*$ represents a semiample Cartier divisor on $\hY$ is immediate from Proposition~\ref{prop:divpoly}, and the formula concerning the coefficients of $\hD$ follows from Proposition~\ref{prop:dcoeff}. That $\hD$ is a polyhedral divisor follows immediately from this representation. Furthermore, a straightforward calculation shows that $\deg \Psi[\hu](u')$ is strictly positive for $\hu$ in the interior of $\hwt$ and $u'$ in the interior of $\Box[\hu]$, so $\hD$ is in fact proper. Finally, by construction,
$$\bigoplus_{u'\in\Box[\hu]}L(\D(u'+s^*(\hu)))=L(\hD(\hu))$$
so $X(\D)=X(\hD)$.
\end{proof}

\begin{ex}[The cone over $G(2,4)$ revisited]
Once again, let $X$ be the cone over $G(2,4)$.
If we consider the upgraded p-divisor of Example \ref{ex:grassupgrade} from Figure \ref{fig:cube}, we can downgrade the torus action of the four-dimensional torus on $X$ to the two-dimensional torus with cocharacter lattice spanned by $e_1$ and $e_2$. This yields exactly the starting data before we performed the upgrade, i.e. the divisorial fan from Figure \ref{fig:grassfan} and downgraded polyhedral divisor from Figure \ref{fig:grasscoeff}.

Note that instead of applying Theorem \ref{thm:downgrade} to calculate the downgraded p-divisor, we could instead have considered $X$ as a hypersurface in $\CC^6$.  By performing a \emph{toric} downgrade on $\CC^6$ as in Section 11 of \cite{altmann:06a}, we get a p-divisor describing the restricted torus action on $\CC^6$. By pulling back this p-divisor to a suitable quotient of $X$, we also get a downgraded p-divisor for $X$. The disadvantage of this approach is that finding the suitable quotient of $X$ requires some calculation, and does not automatically give this quotient the structure of a $T'$-variety.    
\end{ex}

\subsection{Pathologies in the General Case}\label{sec:downgradetwo}
Consider now $Y$ any normal semiprojective variety, and $\D$ a p-divisor with locus $Y$ and coefficients in $N_\QQ$. For a subtorus $\hT\subset T=N\otimes\CC^*$, we would like a normal semiprojective variety $\hY$ and p-divisor $\hD$ on $\hY$ describing $X(\D)$ as a $\hT$-variety. The main difficulty is to find such a variety $\hY$. As an ansatz, one might try to construct $\hY$ as a $T/\hT$-variety from some divisorial fan on $Y$, just as we have done in the complexity-one case. However, as the following extended example shows, in general this cannot work:

\begin{ex}[A downgrade with difficulties]
Consider the lattice $N=\ZZ^2$ with basis $\{e_1,e_2\}$ and corresponding torus $T=\CC^*\otimes N$. This torus splits into a product $T=\hT\times T'$, with $\hT$ and $T'$ corresponding to the first, respectively second, $\ZZ$-summand in $N$. Now let $Y=\AA^2=\spec \CC[x,y]$ and consider the prime divisors $D_x=V(x)$ and $D_y=V(y)$. Since $Y$ is smooth and affine, the polyhedral divisor
$$\D=\conv\{0,e_2\}\otimes D_x+\conv \{0,e_1+e_2\}\otimes D_y$$
is in fact proper; we denote the corresponding variety by $X=X(\D)$. We shall see that it is not possible to describe $X$ as a $\hT$-variety via any p-divisor $\hD$ on any $T'$-variety coming from a divisorial fan on $Y$. 

To begin with, notice that $Y$ is actually toric (with torus $T_Y$), and $\D$ is $T_Y$-invariant. We can thus invoke (a simple form of) Theorem \ref{thm:tD} to get that $X$ is the $\tT=T\times T_Y$-variety corresponding to the cone $\tsigma$ whose rays are generated by the columns of
{\small{$$
\left(
\begin{array}{c c c c}
0&0&0&1\\
0&1&0&1\\
1&1&0&0\\
0&0&1&1\\
\end{array}
\right)
$$}}
in $M\oplus\ZZ^2$. This cone is just a simplex, whose primitive generators have determinant one, so in fact $X=\AA^4$.

Now, although we don't directly know how to downgrade the $T$-action of $X$ to a $\hT$-action, we can use the standard toric downgrade procedure of \cite[Section 11]{altmann:06a} to downgrade the $\tT$-action of $X$ to a $\hT$-action. The Chow quotient of $X$ by $\hT$ is the toric variety $\hY$ corresponding to the image fan $\Sigma_{\hY}$ of $\tsigma$ under the projection which kills $e_1$. If $v_1,\ldots,v_4$ denote the images of the columns of the above matrix, and $v_0=(1,1,1)$, then $\Sigma_{\overline Y}$ has four top-dimensional cones, generated by
\begin{align*}
\{v_0,v_1,v_2\},\qquad \{v_0,v_1,v_3\},\\
\{v_0,v_3,v_4\},\qquad \{v_0,v_2,v_4\},
\end{align*}
see the affine slice of $\Sigma_{\hY}$ pictured in Figure \ref{fig:imagefan}(a).  Note that the ``extra'' ray generated by $v_0$ arises due to common refinement, for example, as the intersection of the cones generated by $\{v_1,v_4\}$ and $\{v_2,v_3\}$.
A p-divisor $\hD$ which describes $X$ as a $\hT$-variety is then 
$$
\hD=\{1\}\otimes D_4+[0,1]\otimes D_5,
$$
where $D_i$ is the toric divisor on $\hY$ corresponding to $v_i$.
\begin{figure}
\subfigure[Slice of $\Sigma_{\hY}$]{\fanslicea}
\subfigure[Slice of $\Sigma_Z$]{\fansliceb}
\caption{Chow quotients $\hY=X\gquot\hT$ and $Z=\hY\gquot T'$}\label{fig:imagefan}
\end{figure}

The torus $T'$ still acts on $\hY$. The Chow quotient $Z$ of $\hY$ by $T'$ is now just the toric variety corresponding to the image fan $\Sigma_Z$ of $\Sigma_Y$ under the projection which kills $e_2$. $\Sigma_Z$ thus has just two top-dimensional cones, generated by $\{(1,0),(1,1)\}$ and $\{(0,1),(1,1)\}$, see Figure \ref{fig:imagefan}(b). Thus, $Z$ is the blowup of $Y=\AA^2$ at the origin. 
In fact, we have the following commutative diagram:
$$
\xymatrix{
 X(\hD) \ar@{-->}[d] \ar@{=}[r] &X(\D)\ar[dd]\\
 \hY\ar@{-->}[d] & \\
  Z\ar[r] & Y. }
$$
$\hY$ corresponds to the divisorial fan $\overline\dfan$ on $Z$ induced by the four p-divisors
\begin{align*}
&[0,1]\otimes D'_x+\emptyset\otimes D'_y+\{1\}\otimes E,\\ 
&\emptyset\otimes D'_x+[0,1]\otimes D'_y+\{1\}\otimes E,\\
&\{1\}\otimes D'_x+\{1\}\otimes D'_y+[1,2]\otimes E,\\
&[0,1]\otimes E,
\end{align*}
where $E$ is the exceptional divisor, and $D'_x$ and $D'_y$ are the strict transforms of $D_x$ and $D_y$.
 
Now, suppose that there is some divisorial fan $\dfan$ on $Y$ and some p-divisor $\E$ on $X(\dfan)$ such that $X=X(\E)$ as a $\hT$-variety. Since $\hD$ is minimal, there must be a dominant $T'$-equivariant map $X(\dfan)\to \hY$ such that $\E$ is the pullback of $\hD$ up to a principal polyhedral divisor, see Theorem \ref{thm:morphisms} and \cite[Theorem 8.8]{altmann:06a}. Furthermore, the map $X(\dfan)\to\hY$ corresponds to a 
triple $(\psi,F,\mathfrak{f})$ as in Definition \ref{def:map}, where in particular $\psi\colon W\to Z$ is some birational map and $W$ is a modification of $Y$. For $\E$ to be the pullback of $\hD$ up to some shift, we must have that the image of the pullback of $\hD(u)$ in $W$ is linearly equivalent to the image of the pullback of $\E(u)$ for all $u\in \ZZ$. Since $Y$ has trivial Picard group, the latter is always trivial.  However, the former is not: the image of $\hD(-1)$ is linearly equivalent to $E$. Thus, $\hD$ cannot pull back to something equivalent to $\E$. This shows that there is no divisorial fan $\dfan$ on $Y$ such that there is a p-divisor $\E$ on $X(\dfan)$ describing $X$ as a $\hT$-variety.
\end{ex}

\section{Upgrades with $\Pic(Y')=\ZZ$}\label{sec:correct}
In this section, we will look at upgrades in an especially simple case.
Consider a polyhedral divisor $\D$ with locus some normal semiprojective variety $Y$ such that $\Pic(Y)=\ZZ$ with $\QQ_{\geq 0}$ the cone of effective divisors. Then one can easily \emph{correct} $\D$ to make it semiample. Indeed, take $\deg \D=\sum \deg(P)\cdot \D_P$, and
$$
\D':=\sum(\D_P+\QQ_{\geq 0} \cdot \deg \D)\otimes P. 
$$ 
Then $\D'$ is a semiample polyhedral divisor, and $X(\D)=X(\D')$. What we have done is simply expand the tailcone of $\D$ so that $\D$ always has nonnegative degree, and is thus always semiample. $\D'$ is in fact a p-divisor if and only if $\deg \D'\neq \tail \D'$ and $\tail \D'$ is pointed.

\begin{rem}
Suppose that instead of $\Pic(Y)=\ZZ$, we just have $\NE(Y)=\ZZ$. Then the corrected p-divisor $\D'$ might not be semiample, but it is nef.
\end{rem}

We can couple the above correction with our upgrade procedure to produce semiample polyhedral divisors, even when the starting polyhedral divisor has no positivity properties. The setup is as follows. Let $\D$ be a polyhedral divisor on some variety $Y$, where $Y=X(\dfan)$ is a $T'$-variety for some divisorial fan on a normal semiprojective variety $Y'$. This is our standard upgrade situation, although we make no assumptions on $\D$ or on $Y$. Instead, we only assume that $\Pic(Y')=\ZZ$ with $\QQ_{\geq 0}$ the cone of effective divisors.
As in Section \ref{sec:upgrade}, we write $\D=\sum \Delta_\rho\otimes D_\rho+\sum \Delta_{P,v}\otimes \mu(v)D_{P,v}$, and set 
$$
{\tsigma} = \pos \big\{ (\sigma \times \{0\}) \cup 
\bigcup_{\rho\in\xray(\dfan)} (\Delta_\rho \times \{v_\rho\}) \big\}
$$
and
$$
\Delta_P = \conv\big\{\Delta_{P,v} \times \{v\} \mid v \in \xvert_P(\dfan) \big\} + \tsigma.
$$
\begin{prop}\label{prop:upcorrect}
Let $\widehat\sigma=\pos \sum_P\deg(P)\cdot\Delta_P$. Then
$$
\widehat{\D}=\sum_P (\Delta_P+\widehat\sigma)\otimes P
$$
is semiample, and $X(\D)=X(\widehat{\D})$. 
\end{prop}
\begin{proof}
The  claims follow directly from the above and Theorem \ref{thm:tD}.
\end{proof}
\begin{rem}
Again, if we just assume that $\NE(Y')=\ZZ$ with $\QQ_{\geq 0}$ the cone of effective divisors, then the above proposition is true if we replace semiample with nef.
\end{rem}

\begin{ex}[A non-contraction-free $\PP^2$ revisited]
Consider again Example \ref{ex:nonsnake}, where we had a non-contraction-free divisorial fan $\dfan$ with $Y=X(\dfan)=\PP^2$, together with a p-divisor $\D$ on $Y$. We had seen that if we upgrade $\D$ to $\tD$, the resulting $\tD$ was not proper, but when we replaced $\dfan$ by a contraction-free $\dfan'$ and pulled by $\D$, the resulting upgraded polyhedral divisor was proper; see Figure \ref{fig:nonsnakeup}.

Instead of replacing $\dfan$ by $\dfan'$, we could also make the correction outlined in this section, replacing $\tsigma$ by $\widehat{\sigma}$. A simple calculation shows that this yields exactly the same p-divisor as when moving from $\dfan$ to $\dfan'$.
\end{ex}

\section{Cox Rings of $T$-Varieties}\label{sec:cox}
Let $X$ be a normal, $\QQ$-factorial complete variety with $\Cl(X)$ finitely generated and torsion free. The $\Cl(X)$-graded group 
$$
\Cox(X):=\bigoplus_{D\in\Cl(X)} L(D)
$$
carries a natural ring structure and is called the \emph{Cox ring} of $X$, see for example \cite{hausen:09a}. If $\Cox(X)$ is a finitely generated $\CC$-algebra, then $X$ is called a \emph{Mori dream space} (MDS) and has many nice geometric properties, see \cite{hu:00a}.

In \cite{altmann:09a}, K. Altmann and J. Wi\'sniewski showed how to describe the spectrum of $\Cox(X)$ in terms of a p-divisor when $X$ is an MDS, although a completely explicit description of this p-divisor is not always immediate. If $X$ is a rational complexity-one $T$-variety, then it is automatically an MDS. In this setting, K. Altmann and L. Petersen described the spectrum of $\Cox(X)$ as a complexity-one $T$-variety via a p-divisor on a curve by considering the action of the Picard torus together with the torus acting on $X$, see \cite{altmann:10a}.\footnote{Altmann and Petersen also include the case where $\Cl(X)$ has torsion, in which case the corresponding p-divisor lives on a covering of $\PP^1$.} In the following, we use our upgrade procedure coupled with the discussion of the previous section  to shed new light on their result and slightly generalize it.

Let $\dfan$ be a divisorial fan in a lattice $N$ on some normal projective variety $Y$ with $\Cl(Y)=\ZZ$ and $\QQ_{\geq 0}$ the effective cone of divisors. Assume that $X=X(\dfan)$ is complete, $\QQ$-factorial, and with $\Cl(X)$ finitely generated and torsion-free. Let $\mcP$ be some nonempty finite set of prime divisors on $Y$ including all $P$ with $\dfan_P$ nontrivial. Let $\V$ be the set of all $(P,v)$ with $P\in\mcP$, $v\in\xvert_P(\dfan)$, and set $\R=\xray(\dfan)$. Then after choosing a section $t^*$ and corresponding cosection $s$ we have the following exact sequence, cf. \cite[Corollary 2.3]{altmann:10a}:
\begin{equation*}
\xymatrix{
0 \ar[r]& \Cl(X)^* \ar[r]  &\ZZ^{(\V\cup\R)} \ar^{\pi}[r] \ar^{s}@/^1pc/[l] & \ZZ^\mcP/\ZZ\oplus N \ar[r] \ar^{t^*}@/^1pc/[l]& 0
}
\end{equation*}
Here,  $e(\cdot)$ are the natural basis elements of $\ZZ^{(\V\cup\R)}$ and $\ZZ^\mcP$,
and 
\begin{align*}
\pi(e(P,v))&=\mu(v)\overline{e(P)}+\mu(v)v;\\
\pi(e(\rho))&=v_\rho.
\end{align*} 
The quotient $\ZZ^\mcP/\ZZ$ is induced by the inclusion $\ZZ\hookrightarrow \ZZ^\mcP$ mapping $1$ to $\sum \deg(P)e(P)$.
A straightforward calculation then shows $\Spec \Cox(X)=X(\D^\cox)$, where
$$
\D^\cox=\sum s(e(P,v))\otimes D_{P,v}+\sum s(e(\rho))\otimes D_\rho
$$
is a polyhedral divisor on $X$. In general, it is very far from being proper.

Nonetheless, we can still apply our upgrading procedure to get a polyhedral divisor $\tD^\cox$ on $Y$:
\begin{align*}
	\tsigma&=\pos \bigcup_{\rho\in\xray(\dfan)} s(e(\rho))\times \{v_\rho\}\\
\tD^\cox&=\sum \tD_P^\cox\otimes P\\
\tD_P^\cox&=\conv \{s(e(P,v)/\mu(v))\times \{v\}\ |\ v\in\xvert_P(\dfan)\}+\tsigma\subset\Cl(X)^*\oplus N.
\end{align*}
Clearly we still have $\Spec \Cox(X)=X(\tD^\cox)$.
We can describe the coefficients of $\tD^\cox$ in a even nicer manner.
Now, we have an inclusion $\Cl(X)^*\times N\hookrightarrow \ZZ^{(\V\cup\R)}$ induced by $t^*$; this gives rise to
\begin{equation*}
\xymatrix{
0 \ar[r]& \Cl(X)^*\oplus N \ar[r]  &\ZZ^{(\V\cup\R)} \ar[r] \ar^{\widetilde s}@/^1pc/[l] & \ZZ^\mcP/\ZZ \ar[r] \ar^{{\widetilde t}^*}@/^1pc/[l]& 0.
}
\end{equation*}
Then
$$
\tD_P^\cox=\widetilde{s}\big(\conv\{e(P,v)/\mu(v)\}\big)+\tsigma=\conv \{e(P,v)/\mu(v)\}+\QQ_{\geq 0}^\R-{\widetilde t}^*(e(P)).
$$ 

This polyhedral divisor is still not proper.  However, since $\Cl(Y)=\ZZ$, we can modify it as in the previous section so that it becomes proper. Indeed, we again take $\deg \tD^\cox=\sum \deg(P)\cdot \tD_P^\cox$, and set $\widehat\sigma=\QQ_{\geq 0}\cdot  \deg \tD^\cox$. We then set 
$$
\widehat{\D}^\cox=\sum (\tD_P^\cox+\widehat\sigma)\otimes P.
$$ 
Bigness is easy to check, and one can conclude that $\widehat{\D}^\cox$ is a p-divisor with $\Spec \Cox(X)=X(\widehat{\D}^\cox)$.  A quick comparison shows that for $Y$ a curve, this is exactly the p-divisor of \cite[Theorem 1.2]{altmann:10a}. However, we are no longer restricted to calculating the p-divisor of Cox rings for rational complexity-one $T$-varieties. For example, we can compute the p-divisor of the Cox ring of any $T$-variety $X(\dfan)$, where $\dfan$ is a divisorial fan on some weighted projective space or Grassmannian.

\section{Upgrading deformations of toric varieties}\label{sec:def}

Let $X$ be a toric variety with embedded torus $\tT$,
given by a polyhedral cone $\delta$ in a lattice $\tN$.
For a degree $r \neq 0$ in the dual lattice $\tM$, 
K. Altmann has shown how certain Minkowski decompositions
$$
\delta_r := \delta \cap [r=1] = \Delta_0 + \dotsm + \Delta_l
$$
correspond to $l$-parameter deformations
$\cX \to \AA^l$ of $X$
with deformation parameters of degree $r$~\cite{altmann:00a}.
For degrees $r \in \delta^\vee$, $\cX$ itself is a toric variety,
while the construction is more complicated for general $r$.

Let $r_0 \in \tM$ be primitive such that $r = kr_0$, $k \in \NN$. 
Then the deformation $\cX$ is invariant with respect to the
torus $T := \ker r_0 \subset \tT$, i.e. $T$ acts on each fiber.
To simplify notation, we choose a splitting $\tT = T \times T'$,
so that $r_0$ is the projection onto $T' = \CC^*$, and $\tN = N 
\oplus \ZZ$. The dual of the weight cone of the $T$-action 
on $X$ is $\sigma := \delta \cap (N_\QQ,0)$, and with
\begin{align*}
  (\Delta^+,1) & := \delta_{r_0} = \delta \cap (N_\QQ,1) \\
  (\Delta^-,-1) & := \delta_{-r_0} = \delta \cap (N_\QQ,-1),
\end{align*}
$X$ has the p-divisor
$\D = \Delta^+ \otimes 0 + \Delta^- \otimes \infty$
on $\PP^1$. Note that $\Delta^-$ is the empty set in the case
$r \in \delta^\vee$.
In an earlier article, the present authors have shown how to express
the family $\cX$ as a $T$-variety, as summarized 
below~\cite{ilten:09b}.

\begin{definition}
An \emph{admissible} Minkowski decomposition of $\delta_r$ consists
of polyhedra $\Delta_0, \dotsc, \Delta_l$ in $N_\QQ$ with 
common tailcone $\sigma$ such that
$$
\delta_r = (\Delta_0,\sfrac{k}) + (\Delta_1,0) + \dotsm + (\Delta_l,0)
$$
and that for each $u \in M \cap \sigma^\vee$, at most one of the 
$u$-faces of the summands does not contain lattice points.
\end{definition}

\begin{rem}
For non-primitive $r$, that is if $k > 1$, 
$(\Delta_0,\frac{1}{k})$ doesn't contain lattice points, so the
polyhedra $\Delta_1,\dotsc,\Delta_l$ must be lattice polyhedra.
Note that this decomposition induces a decomposition
$$
\delta_{r_0} = (k\Delta_0,1) + (k\Delta_1,0) + \dotsm + (k\Delta_l,0)
$$
in the primitive degree $r_0$.
\end{rem}

\begin{prop}[{\cite[Section 3]{ilten:09b}}]\label{prop:pdivdef}
Let $\cX \to \AA^l$ be the deformation of $X$ associated to an 
admissible decomposition of $\delta_r$.
Let $Y = \PP^1 \times \AA^l$ and consider the divisors
$D_0 = V(x_1)$, $D_i = V(x_1^k - y_i x_0^k)$ for $1 \le i \le l$,
and $D_\infty = V(x_0)$.
Then $\cX$ is given by the p-divisor
$$
\E = k\Delta_0 \otimes D_0 + \Delta_1 \otimes D_1 + \dotsm + 
\Delta_l \otimes D_l + \Delta^- \otimes D_\infty
$$
on $Y$. The structure map of the family corresponds to the
projection $Y \to \AA^l$.
\end{prop}

\begin{figure}
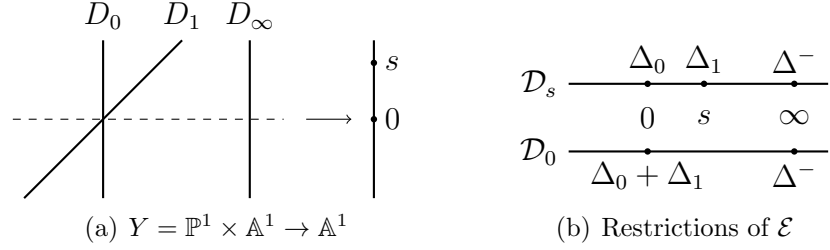

\centering
\subfigure[$Y = \PP^1 \times \AA^1 \to \AA^1$]{\tvardefBase}
\hspace*{5ex}
\subfigure[Restrictions of $\E$]{\tvardefDivisors}
\caption{One-parameter toric deformation as deformation of 
p-divisors.}
\label{fig:tvardef}
\end{figure}

The situation is summarized in Figure~\ref{fig:tvardef} for the
case $r=1$.
We now show how to upgrade the $T$-action on $\cX$ to a $\tT$-action:

\begin{theorem}\label{thm:upgradedef}
Let $\cX \to \AA^l$ be the deformation of $X$ associated with an 
admissible decomposition of $\delta_r$. Then the $\tT$-action
on $X$ extends to $\cX$, with weight cone dual to
$$\tsigma = \delta \cap [r_0 \ge 0].$$
To describe $\cX$ as a $\tT$-variety, consider
$\PP^l = \Proj \CC[y_0, \dotsc, y_l]$, and let $Z$ be the 
blowup of $\PP^l$ at the point $O = (1:0:\dotsm:0)$, with
the divisors $P_0 = V(y_0)$, $P_i = V(y_0 - y_i)$ and the
exceptional divisor $Q$. Then $\cX$ is described by
the p-divisor on $Z$ with tailcone $\tsigma$ and coefficients
\begin{align*}
  \Delta_{P_0} &= (\Delta_0,\sfrac{k}) + \tsigma &
  \Delta_{P_i} &= (\Delta_i,0) + \tsigma &
  \Delta_Q     &= \conv\{(\sfrac{k}\Delta^-,-\sfrac{k}), (0,0)\}
                  + \tsigma.
\end{align*}
\end{theorem}
Note that if $l = 1$, we have $Z = \PP^1$ and $Q = O = V(y_1)$.

\begin{proof}
By Proposition~\ref{prop:pdivdef}, we can describe $\cX$ by the
p-divisor
$$
\E = k\Delta_0 \otimes D_0 + \Delta_1 \otimes D_1 + \dotsm + 
\Delta_l \otimes D_l + \Delta^- \otimes D_\infty
$$
on $Y = \PP^1 \times \AA^l$.
In order to upgrade the $T$-action, we will first
need to express $\E$ as an invariant p-divisor with respect
to an appropriate action of $T'$ on $Y$.

Since $\cX \to \AA^l$ should be $\tT$-equivariant if $\tT$ acts
on $\AA^l$ with weights $r$, it follows that
$\deg y_i = r \in \tM$. For $D_i = V((\frac{x_1}{x_0})^k - y_i)$
to be $T'$-invariant, we want $\deg \frac{x_1}{x_0} = r_0$.
Thus, we consider the $T'$-action on $Y$ given by the weights $(1,k,\dotsc,k)$.

Now by Proposition~\ref{prop:upgradedef} below, $\E$ pulls back to
a $T'$-invariant p-divisor $\E'$ on a contraction-free $T'$-variety
$\tY = X(\dfan)$ over $Z$, where $\dfan$ has slices
\begin{align*}
\dfan_{P_0} &= [\sfrac{k},\infty) &
\dfan_{P_i} &= [0,\infty) &
\dfan_{Q}   &= [-\sfrac{k},0] \cup [0,\infty).
\end{align*}
The positive half-line $\rho$ is the only ray in $\tail \dfan$, and
expressing $\E'$ in the form required by Theorem~\ref{thm:tD}, we get
\begin{align*}
\E' &= \Delta^+ \otimes D_\rho \\
    &+ \Delta_0 \otimes kD_{(P_0,\frac{1}{k})} \\
    &+ \Delta_1 \otimes D_{(P_1,0)} + \dotsc
       + \Delta_l \otimes D_{(P_l,0)} \\
    &+ \sfrac{k}\Delta^- \otimes kD_{(Q,-\frac{1}{k})}.
\end{align*}
An application of Theorem~\ref{thm:tD} completes the proof.
\end{proof}

\begin{prop}\label{prop:upgradedef}
Let $\E$ on $Y = \PP^1 \times \AA^l$ be as in the proof of
Theorem~\ref{thm:upgradedef}. 
Recall that $Y$ has coordinates $x_0, x_1$ and $y_1, \dotsc, y_l$.
Let $Z = \Bl_0\PP^l$ be the blowup of $\PP^l$ at the origin $O$,
where $\PP^l$ has homogeneous coordinates $y_0 = (\frac{x_1}{x_0})^k$
and $y_i$, $1 \le i \le l$.
\begin{enumerate}
\item
Denote by $P_0 = V(y_0)$ the hyperplane at infinity, and
let $Q$ be the exceptional divisor.
Then as a $T'$-variety, $Y$ is given by the divisorial
fan $\dfan$ on $Z$ generated by
\begin{align*}
\D^+ &= [\sfrac{k},\infty) \otimes P_0 \\
\D^- &= \emptyset \otimes P_0 + [-\sfrac{k},0] \otimes Q.
\end{align*}

\item
Let $P_i = V(y_0 - y_i)$ be the image of $D_i$ in $Z$ for $i \ge 1$.
The prime divisors $D_0$, $D_i$ and $D_\infty$
on $Y$ are $T'$-invariant with
\begin{align*}
  D_0 &= D_{(P_0,\frac{1}{k})} &
  D_i &= D_{(P_i,0)} &
  D_\infty &= D_{(Q,-\frac{1}{k})}.
\end{align*}
In particular, $\E$ is a $T'$-invariant p-divisor.
\item
$Y = \PP^1 \times \AA^l$ becomes contraction-free over $Z$ by 
blowing up to $\tY$ at the origin and along
(the strict transform of)
$\PP^1 \times \{0\}$. This introduces two
(exceptional) divisors $E = D_\rho$, $\rho = \QQ_{\ge 0}$ and
$D_{(Q,0)}$.
\newcommand{\td}{\widetilde{D}}
The divisors $D_i$ pull back to 
\begin{align*}
  \td_0 &= D_{(P_0,\frac{1}{k})} + E &
  \td_i &= D_{(P_i,0)} + E &
  \td_\infty &= D_{(Q,-\frac{1}{k})}.
\end{align*}
\end{enumerate}
\end{prop}

\begin{proof}
The first statement is a straightforward application of the toric downgrade
procedure. The second follows directly from the characterization of Weil
divisors on $T$-varieties. For the last claim, we just need to check how
the divisors $D_i$ intersect the centers of the blowups. The origin
of $\PP^l$ is contained in all $D_i$ but not in $D_\infty$, while 
none of the interesting divisors contain $\PP^1 \times \{0\}$.
\end{proof}

\begin{rem}
The coefficients of the $P_i$ provide a Minkowski decomposition
$$
\Delta_{P_0} + \dotsm + \Delta_{P_r} = \delta \cap [r \ge 1],
$$
while $\Delta_Q = \delta \cap [r \ge -1]$.
\end{rem}

\begin{rem}
To describe the $\tT$-equivariant structure map $\cX \to \AA^l$,
let $[1,\infty) \otimes H$ on $\PP^{l-1}$ be a p-divisor for $\AA^l$
as a $\CC^*$-variety, where $H$ is any hyperplane in $\PP^{l-1}$.
Then if $\pi \colon Z \to \PP^{l-1}$ is the
projection from the (blown-up) origin $O$, the structure map
corresponds to the triple $(\pi,r,\mathfrak{f})$, where
$\mathfrak{f}(1)$ is the principal divisor $P_0 - \pi^*H - Q$
on $Z$.
\end{rem}

\begin{ex}[Toric total spaces]
If $\Delta^-$ is empty, the coefficient of $Q$ is trivial,
and the upgraded divisor may be linearly translated to one supported
in the coordinate hyperplanes $P'_0,\dotsc,P'_r$. This is a configuration
of torus invariant divisors on the toric variety $Z$, which shows that
in this case $\cX$ is toric itself.
\end{ex}

\begin{figure}
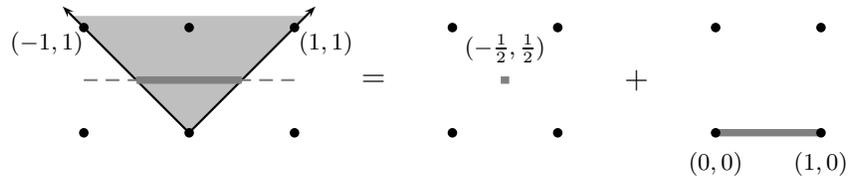

\minkdecomp
\caption{A Minkowski decomposition for the quadric cone}\label{fig:minkdecomp}
\end{figure}

\begin{figure}
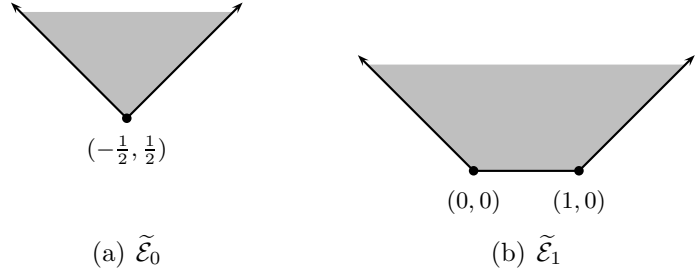

\subfigure[$\widetilde{\E}_0$]{\minkda}
\subfigure[$\widetilde{\E}_1$]{\minkdb}
\caption{Upgraded coefficients for the total space $\cX$}\label{fig:minkdecompup}
\end{figure}

\begin{ex}[An $A_1$-singularity]
Consider the quadric cone $X = V(uw-v^2) \subset \AA^3$. As a toric
variety, this is given by $\delta = \pos \{(1,1), (-1,1)\}$, and the
variables $u,v,w$ have weights $[1,1], [0,1], [-1,1]$. We can
deform $X$ in degree $r = [0,2]$ to $\cX = V(uw-v^2-y_1) \cong 
\AA^3$, corresponding to the Minkowski decomposition
$$
\delta_r = \left([-\sfrac{2},\sfrac{2}] , \sfrac{2}\right)
         = (\Delta_0,\sfrac{2}) + (\Delta_1 ,0)
$$
with $\Delta_0=\{-1/2\}$ and $\Delta_1=[0,1]$, see Figure \ref{fig:minkdecomp}.
Since $\Delta^- = \emptyset$, we may take $Y = \AA^1 \times \AA^1$
with coordinates $x$ and $y$.
Then as a $T$-variety, $\cX$ is given by the p-divisor $\E$ on $Y$,
where $P_0 = V(x)$ has coefficient $\{-1\}$, and
$P_1 = V(y - x^2)$ has coefficient $[0,1]$.
Upgrading this divisor gives
$$
\widetilde{\E}=\big((\Delta_0,\sfrac{2}) + \delta\big) \otimes \{0\}
    + \big((\Delta_1,0) + \delta\big) \otimes \{1\}
$$
on $\PP^1$, see Figure \ref{fig:minkdecompup}.
\end{ex}

\begin{rem}
The same approach allows for the description of certain mixed deformations,
where different multiples of the same primitive degree $r_0$ occur.
Namely, consider an admissible decomposition
$$
\Delta^+ = \Delta_0 + k_1\Delta_1 + \dotsm + k_l\Delta_l
$$
where the multiplicities $k_i$ are part of the data, i.e. 
$\Delta_0 + 2\Delta_1$ differs from $\Delta_0 + \Delta_1 + 
\Delta_1$. This determines a $T$-invariant
deformation of $X$~\cite{ilten:09b}, given by the $T'$-invariant 
divisor
$$
\E = \Delta_0 \otimes D_0 + \Delta_1 \otimes D_1 + \dotsc + 
\Delta_l \otimes D_l + \Delta^- \otimes D_\infty
$$
on $Y = \PP^1 \times \AA^l$, where $D_0$ and $D_\infty$ are as 
above and $D_i = V(x_1^{k_i} - y_ix_0^{k_i})$ are adapted to
varying multiplicities. $T'$ acts on $Y$ with weights 
$(1,k_1,\dotsc,k_l)$.

If $k$ is the greatest common divisor of the $k_i$, and $k_i = h_ik$,
then we can upgrade $\E$ to a p-divisor on the blowup of
the weighted projective space $\PP(1,h_1,\dotsc,h_l)$ at the
origin. With $P_i = V(y_0^{h_i} - y_i)$, this divisor has
coefficients
\begin{align*}
  \Delta_{P_0} &= (\Delta_0,\sfrac{k}) + \tsigma &
  \Delta_{P_i} &= (\Delta_i,0) + \tsigma &
  \Delta_Q     &= \conv\{(\sfrac{k}\Delta^-,-\sfrac{k}), (0,0)\}
                  + \tsigma.
\end{align*}
\end{rem}

\bibliographystyle{alpha}
\bibliography{upgrade}

\end{document}